\documentstyle[12pt]{article} 
\makeatletter
\def\mineappendix{
        \setcounter{section}{1}
        \setcounter{subsection}{0}
        \def\thesection{\Alph{section}}
        \def\sectionap{\@startsection  {section}{1}{\z@}
                        {-3.5ex plus-1ex minus-.2ex} {0ex plus.2ex}
                        {\reset@font\Large\bf  Appendix:  \, }
                        } 
        }
\makeatother
\def\Proclaim #1. #2\par{\bigbreak\noindent{\sc#1.\enspace}{\it#2}\par}
\newtheorem{lem}{Lemma}[section]
\newtheorem{cor}[lem]{Corollary} 
\newtheorem{thm}[lem]{Theorem} 
\newtheorem{pro}[lem]{Proposition} 
\newtheorem{exa}[lem]{Example} 
\newtheorem{defi}[lem]{Definition}  
\newtheorem{rem}[lem]{Remark}

\title{Elliptic Gromov-Witten Invariants And Virasoro Conjecture}
\author{Xiaobo Liu}
\date{}

\begin{document}
\maketitle

\input amssym.def

The Virasoro conjecture predicts that the generating function
of Gromov-Witten invariants is annihilated by infinitely many
differential operators which form a half branch of the Virasoro
algebra. This conjecture was proposed by Eguchi, Hori and Xiong \cite{EHX2}
and also by S. Katz \cite{Ka} (see also \cite{EJX}).
It provides a powerful tool in the computation of Gromov-Witten invariants.
In \cite{LT}, the author and Tian proved the genus-0 part of the
Virasoro conjecture.
The main purpose of this paper is to study the genus-1 part of this 
conjecture. 

The system of Gromov-Witten invariants relevant to this paper are the
so called descendant Gromov-Witten invariants. These invariants
 arose in the theory of topological 
sigma model coupled to gravity \cite{W2}. Mathematical definition for
such invariants was given in \cite{RT2} for semipositive symplectic
manifolds. Using virtual moduli cycles (\cite{LiT1}, \cite{LiT2}, and
also \cite{BF}),
these invariants can also be defined for all compact symplectic
manifolds and, in purely algebraic geometric setting, for smooth projective
varieties. In this paper, we consider descendant Gromov-Witten invariants
for a smooth projective variety $V$. For simplicity, we assume that
$H^{\rm odd}(V, {\Bbb C}) = 0$. Fix a basis $\{\gamma_{1}, \ldots,
\gamma_{N}\}$ of $H^{*}(V, {\Bbb C})$ with $\gamma_{1}$ equal to the
identity of the cohomology ring of $V$ and 
$\gamma_{\alpha} \in H^{p_{\alpha}, q_{\alpha}}(V, {\Bbb C})$ for every
$\alpha$. For any non-negative integer $g$ and 
$A \in H_{2}(V, {\Bbb Z})$, let 
$\left< \tau_{n_{1}, \alpha_{1}} \ldots \tau_{n_{k}, \alpha_{k}}
		\right>_{g, A}$
be the genus $g$ degree $A$ 
descendant Gromov-Witten invariants associated with cohomology
classes $\gamma_{\alpha_{1}}, \ldots, \gamma_{\alpha_{k}}$ and non-negative
integers $n_{1}, \ldots, n_{k}$ 
(See Section~\ref{sec:defGW} for the definition of Gromov-Witten invariants). 
Summing up the Gromov-Witten invariants over all degrees, we obtain a 
quantity which is called the $k$-point correlators in the theory of
topological sigma model:
\[ \left< \tau_{n_{1}, \alpha_{1}} \ldots \tau_{n_{k}, \alpha_{k}} \right>_{g}
	 \, \, \, := \sum_{A \in H_{2}(V, {\Bbb Z})} q^{A}
	 \left< \tau_{n_{1}, \alpha_{1}} \ldots \tau_{n_{k}, \alpha_{k}}
			\right>_{g, A},\]
where $q^{A}$ belongs to the Novikov ring (i.e. the completion of the
 multiplicative ring
generated by monomials $q^{A} = d_{1}^{a_{1}} \cdots d_{r}^{a_{r}}$
over the ring of rational numbers, where
$\{d_{1}, \cdots, d_{r}\}$ is a fixed basis of $H_{2}(V, {\Bbb Z})$ and
$A = \sum_{i=1}^{r} a_{i}d_{i}$).
The generating function of genus-$g$ Gromov-Witten invariants is defined
by
\[ F_{g}(T) := \sum_{k \geq 0} \frac{1}{k!}
	 \sum_{ \begin{array}{c} {\scriptstyle \alpha_{1}, \ldots, \alpha_{k}}
					 \\
		                {\scriptstyle  n_{1}, \ldots, n_{k}}
		\end{array}}
		t^{\alpha_{1}}_{n_{1}} \cdots t^{\alpha_{k}}_{n_{k}}
          \left< \tau_{n_{1}, \alpha_{1}} \cdots \tau_{n_{k}, \alpha_{k}}
          \right>_{g}, \]
where $T = \{t^{\alpha}_{n} \mid n \in {\Bbb Z}_{+}, \alpha = 1, \cdots, N\}$
is an infinite set of parameters. The space of all parameters $T$
is called the {\it big phase space}. This is an infinite dimensional
space with coordinates $\{t_{n}^{\alpha} \}$. The finite dimensional  subspace 
$\{T \mid t^{\alpha}_{n} = 0 \, \, \, {\rm if} \, \, \,
        n > 0\}$ 
is called the {\it small phase space}. The function $F_{g}$ is understood
as a formal power series of $T$. The generating function for Gromov-Witten
invariants of all genera is defined to be
 \[ Z(T; \lambda) := \exp \sum_{g \geq 0} \lambda^{2g-2} F_{g}(T), \]
where $\lambda$ is another parameter which is used to separate information
from different genera. In topological sigma model, $F_{g}$ is called the
genus-$g$ {\it free energy function} and $Z$ is called the {\it partition
function}.

In \cite{EHX2}, Eguchi, Hori, and Xiong 
constructed a sequence of
linear differential operators, denoted by $L_{n}$ with $n \in {\Bbb Z}$, 
on the big phase space (See Section~\ref{sec:VirOp}).
They  checked that these operators define a representation of the Virasoro
algebra with the central charge equal to the Euler characteristic number
of $V$ under a condition which is equivalent to the vanishing of the Hodge
number $h^{p, q}(V)$ for $p \neq q$. These operators were modified by
S. Katz so that the last condition is not needed. They conjectured that
$L_{n} Z \equiv 0$ for all $n \geq -1$. This conjecture is  called
the {\it Virasoro conjecture} and the equation $L_{n} Z = 0$ is
called the $L_{n}$ constraint. The $L_{-1}$ constraint is the string equation
(cf.\cite{W2}). The $L_{0}$ constraint was discovered by Hori \cite{H}.
Both of these two constraints hold for all manifolds. 
When the underlying manifold is a point, the Virasoro conjecture is equivalent
to a conjecture by Witten \cite{W2} which predicted that the corresponding
 generating function is a $\tau$-function of the KdV hierarchy. Witten's
conjecture was proved by Kontsevich \cite{Kon} and also by Witten \cite{W3}.
For arbitrary manifolds,  we can write 
\[ L_{n}Z(T; \lambda) = 
   \left\{\sum_{g \geq 0} \Psi_{g, n} \lambda^{2g-2} \right\} Z(T; \lambda).\]
The $L_{n}$ constraint is equivalent to $\Psi_{g, n} = 0$ for all
$g \geq 0$. The equation $\Psi_{g, n} = 0$ is called 
{\it genus-$g$ $L_{n}$-constraint}. It is, in general, a non-linear partial
differential equation involving 
all free energy functions $F_{g^{'}}$ with $0 \leq g^{'} \leq g$.
The {\it genus-$g$ Virasoro conjecture} predicts that for
all $n \geq -1$, the genus-$g$ $L_{n}$ constraint is true. 
The genus-0 Virasoro conjecture
was first proved in \cite{LT}. Later, alternative proofs
 were given in \cite{DZ2} and \cite{G2}. 
We will give a brief review to the genus-0 case in Section~\ref{sec:g0Vir}.
The genus-1 Virasoro conjecture
for manifolds with semisimple quantum cohomology was proved in \cite{DZ2}.
There is a discussion of Virasoro conjecture for degree 0 Gromov-Witten
invariants in \cite{GP}.
In this paper, we study the genus-1 Virasoro conjecture without
assuming semisimplicity. The generating functions relevant to the
genus-1 case are $F_{0}$ and $F_{1}$. As in the genus-0 case,
most of our discussions only use basic properties of quantum cohomology,
therefore could be extended to the setting of abstract Frobenius
manifolds.

In most part of this paper (except Section~\ref{sec:g1VirIff}) we
 will deal with the small phase space which can be identified with
$H^{*}(V, {\Bbb C})$. We will write the coordinates 
$t^{\alpha}_{0}$ simply as $t^{\alpha}$ and identify the coordinate
vector fields $\frac{\partial}{\partial t^{\alpha}}$ with cohomology
classes $\gamma_{\alpha}$. The restriction of $F_{g}$ to the
small phase space is denoted by $F_{g}^{s}$. The third derivatives
of $F_{0}^{s}$ defines a ring structure, called the {\it quantum cohomology
ring}, on each tangent space of the small phase space. This enables us
to take product, called the {\it quantum product}, of two vector
fields on the small phase space. There are two special vector
fields on the small phase space. One is $\gamma_{1}$, which is also
the identity element with respect to the quantum product. Another one
is the so called {\it Euler vector field}, which is defined by
\[ E := c_{1}(V) + \sum_{\alpha} (b_{1} + 1 - b_{\alpha}) t^{\alpha} 
		\gamma_{\alpha}, \]
where $c_{1}(V)$ is the first Chern class of $V$ and
 $b_{\alpha}= p_{\alpha} - \frac{1}{2}({\rm dim}_{\Bbb C}V-1)$.
Note that usually the holomorphic dimension $p_{\alpha}$ is replaced by a half
of the real dimension of $\gamma_{\alpha}$. This modification is due
to S. Katz. Let $E^{k}$ be the $k$-th quantum power of $E$, i.e.
\[ E^{k} := \underbrace{E \bullet \cdots \bullet E}_{k}, \]
where $\bullet$ denotes the quantum product. Here we use the convention that
$E^{0} = \gamma_{1}$ and $E^{1} = E$.
 It is perhaps well known that $\{E^{k} \mid k \geq 0\}$ form a half branch
of the Virasoro algebra, i.e.
\begin{equation} \label{eqn:VirEuler}
	[E^{k}, \, E^{m}] = (m-k) E^{m+k-1}. 
\end{equation}
This fact was used in \cite{DZ2} without giving a proof. A proof of this
can be found in \cite{HM} 
(see also the remark after equation (\ref{eqn:derPowerEuler})). 
It is also well known that
$E^{0}F_{1}^{s} = 0$ and $E F_{1}^{s} = {\rm const}$.
In Section~\ref{sec:g0g1}, we will prove the following
\begin{thm} \label{thm:EulerF1}
	For any manifold $V$ and $k > 0$, the genus-1 data
	$E^{k}F_{1}^{s} - (k/2) E^{k-1}E^{2}F_{1}^{s}$
	can be represented by derivatives of $F_{0}^{s}$.
\end{thm}
See Theorem~\ref{thm:g0g1G0} for a more explicit form of this theorem.
According to this theorem, if we know that $E^{2}F_{1}^{s}$ can be represented
by genus-0 data, so does $E^{k}F_{1}^{s}$ for all $k \geq 0$.
When restricted to the small phase space, the genus-1 $L_{1}$ constraint
is equivalent to say that $E^{2}F_{1}^{s}$ is equal to the following
function
\begin{equation} \label{eqn:phi2}
 \phi_{2} := 
	\sum_{\alpha, \beta} \eta^{\alpha, \beta}
	\left\{ - \frac{1}{24} 
	   \nabla^{2}_{E, \, E} (\gamma_{\alpha} \gamma_{\beta} F_{0}^{s})
	+ \frac{1}{2} \left(b_{\alpha} b_{\beta} - \frac{b_{1}+1}{6}\right)
		\gamma_{\alpha} \gamma_{\beta} F_{0}^{s} \right\},
\end{equation}
where $(\eta^{\alpha, \beta})$ is the inverse matrix of the intersection
form on $V$, $\nabla$ is the flat connection and 
$\nabla^{2}_{u, \, v} = \nabla_{u}\nabla_{v} - \nabla_{\nabla_{u}v}$ is
the second covariant derivative.
In Section~\ref{sec:g1VirIff}, we will prove that this last condition implies the genus-1 $L_{n}$ constraints
for all $n \geq 1$.
\begin{thm} \label{thm:phi2g1Vir}
	For any manifold $V$, the genus-1 Virasoro conjecture holds if and 
	only if $E^{2}F_{1}^{s} = \phi_{2}$.
\end{thm}

Theorem~\ref{thm:EulerF1} gives an expression for $E^{k}F_{1}^{s}$ 
in terms of $E^{2}F_{1}^{s}$ and genus-0 data. If we replace
$E^{2}F_{1}^{s}$ in this expression by $\phi_{2}$ and denote the resulting
expression by $\phi_{k}$ for $k > 2$, then $\phi_{k}$ is a function
only involves genus-0 data (see formula (\ref{eqn:phi-k}) and 
Theorem~\ref{thm:phik} 
 for  explicit forms of this function).
 We also define $\phi_{0} = 0$ and $\phi_{1}$
 equal to the constant in $EF_{1}^{s}$. 
Because of the relation~(\ref{eqn:VirEuler}), a necessary condition for
$E^{2}F_{1}^{s} = \phi_{2}$ is the following
\begin{equation} \label{eqn:phiVir}
 E^{k} \phi_{m} - E^{m} \phi_{k} = (m-k) \phi_{k+m-1},
\end{equation}
for all $m$, $k \geq 0$.
In Section~\ref{sec:phiVir}, we will 
prove that this condition is always satisfied.
\begin{thm} \label{thm:phiVir} 
	For any manifold $V$, equality~(\ref{eqn:phiVir}) always holds.
\end{thm}

At each point of the small phase space, the quantum powers of the Euler
vector field span a subspace of the tangent space. The dimension of this
subspace may vary as the base point changes. In an open subset,
this dimension is constant and the quantum powers of the Euler vector
field define an integrable
distribution. Therefore one can talk about leaves of this distribution.
In fact, in a proper sense,
 leaves of any collection of vector fields on a finite dimensional
manifold  are always well defined and are immersed
submanifolds (see \cite{Sussmann}). Each leaf of $\{E^{k} \mid k \geq 0\}$ 
is a finite dimensional smooth submanifold, which may not be flat with respect
to the intersection form, and therefore may not be a Frobenius manifold itself.
On each leaf (restricting to an open subset if necessary), there exists a 
finite number n such that 
\begin{equation} \label{eqn:wrapEuler}
 E^{n+1} = \sum_{k=0}^{n} f_{k} E^{k},
\end{equation}
where $f_{k}$'s are smooth functions on the leaf. Another necessary condition
for $E^{2}F_{1}^{s} = \phi_{2}$ is that, on each leaf, 
\begin{equation} \label{eqn:philinear}
 \phi_{n+1} = \sum_{k=0}^{n} f_{k} \phi_{k}.
\end{equation}
We conjecture that this condition is always satisfied. 
This can be verified easily for manifolds with semisimple quantum cohomology.
In fact, equations (\ref{eqn:phiVir}) and (\ref{eqn:philinear}) are equivalent
to the existence of a local potential function whose derivative along
$E^{k}$ is $\phi_{k}$ for all $k$. For manifolds with semisimple
quantum cohomology, such a potential function exists globally
and can be explicitly expressed in terms of the $\tau$-function of
the isomonodromy deformation (c.f. \cite{DZ2} proof of Proposition 4).

\begin{defi} \label{def:nondeg}
	We say that a manifold $V$ has {\bf non-degenerate} quantum cohomology
if at generic points of the small phase space, there exists an integer
$m \geq 1$ such that $E^{m}$ 
is contained in the linear span of $\{E^{0},  Z_{k} \mid k \geq 0 \}$, where
$Z_{k} := \sum_{i=0}^{n} (E^{k}f_{i})E^{i}$. 
\end{defi}
We have the following 
\begin{thm} \label{thm:nondegVir}
	For any manifold $V$ with non-degenerate quantum cohomology, if
equality (\ref{eqn:philinear}) is satisfied , then
the genus-1 Virasoro conjecture holds.
\end{thm}
This theorem will be proved in Section~\ref{sec:SuffVir}.  
Because of this theorem,  it would be interesting to know which 
manifolds have  non-degenerate quantum cohomology.
We first note that vector fields $ Z_{k}$ 
can also be defined by $f_{i}$'s without taking derivatives 
(see Remark~\ref{rem:fDer} and Remark~\ref{rem:ZkZ0}).
Therefore the property of being non-degenerate can  be checked pointwise.
In Section~\ref{sec:SuffVir}
 we will give some sufficient conditions for the non-degeneracy.
In particular, if the quantum cohomology of a manifold is semisimple,
it must also be
non-degenerate. As a corollary of Theorem~\ref{thm:nondegVir},  
the genus-1 Virasoro conjecture holds for
manifolds with semisimple quantum cohomology. 
This fact was proved before in \cite{DZ2}.
In the approach of \cite{DZ2}, the assumption of semisimplicity was needed
from the very beginning since the canonical coordinates are used throughout
all calculations. While in our approach, this is a corollary of a more 
general result and our assumption of non-degeneracy only comes at the 
last step.
We would also like to make a comparison between 
the non-degeneracy condition and
the semisimplicity condition. If a manifold has semisimple quantum cohomology,
 then at generic points, the powers of the Euler vector field span the 
entire tangent spaces. But for manifolds with non-degenerate quantum
cohomology, this may not
be the case. Even if we assume that the powers of the Euler vector field
span tangent spaces, the non-degeneracy condition is still weaker
than semisimplicity (we will see this through
Lemma~\ref{lem:NonDegSS}, its corollaries, and examples at the end
of Section~\ref{sec:SuffVir}).
Moreover, to verify semisimplicity, we need to know the quantum product of
the Euler vector fields with tangent vectors in all directions. But, to
verify non-degeneracy, we only need to know the quantum powers of the
Euler vector field. Therefore it might be much easier to give a more geometric
characterization. We recall a conjecture by Tian
which predicts that all Fano varieties (which by definition
 have positive first Chern classes) have semisimple quantum 
cohomology \cite{T}. This conjecture was verified for Grassmannians
and complete intersections of low degrees (see \cite{TX}). In general,
it is still an open conjecture. A weaker version of this conjecture
would be that all Fano varieties have non-degenerate quantum cohomology.
 Since the definition of the Euler vector field explicitly
involves the first Chern class, it might be easier to verify this
weaker version of Tian's conjecture.
We would like to study this in another paper.

The author would like to thank V. Kac, G. Tian, and E. Witten for very
helpful discussions. He is grateful to G. Tian for encouragement
during this work and collaboration in the previous work. 
The author is partially
supported by an NSF postdoctoral fellowship.

\section{Preliminaries} \label{sect:prelim}

In this section we recall the definition of Gromov-Witten
invariants, Quantum cohomology, and some well known facts. 
We will also set up notation conventions used in this paper
and define the Virasoro operators.
In Section~\ref{sec:g0Vir}, we give a brief review of the genus-0
Virasoro conjecture.

\subsection{Gromov-Witten invariants}
\label{sec:defGW}
Gromov-Witten invariants are defined via the intersection theory
of moduli spaces of stable maps from Riemann surfaces to a fixed manifold
$V$. For any element $A \in H_{2}(V, {\Bbb Z})$ and non-negative integers
$g$ and $k$,  the moduli space $\overline{\cal M}_{g, k}(V, A)$ is defined
to be the collection of all data
$(C; x_{1}, \ldots, x_{k}; f)$ 
where $C$ is a genus-$g$ projective connected
curve over ${\Bbb C}$ whose only possible singularities are simple
double points, $x_{1}, \ldots, x_{k}$ are smooth points on $C$
(called marked points), and $f$ is an algebraic map from
$C$ to $V$ which is stable with respect to $(C; x_{1}, \ldots, x_{k})$,
(i.e. there is no infinitesimal deformation for this data).
Each marked point $x_{i}$ defines a map, called the $i$-th evaluation
map,
\[ \begin{array}{rccc}
	  ev_{i}: & \overline{\cal M}_{g, k}(V, A) & \longrightarrow & V \\
		& (C; x_{1}, \ldots, x_{k}; f) & \longmapsto & f(x_{i}).
	\end{array} \]
It also defines a line bundle over $\overline{\cal M}_{g, k}(V, A)$,
denoted by $E_{i}$, whose fiber over $(C; x_{1}, \ldots, x_{k}; f)$ 
is $T^{*}_{x_{i}}C$. For any cohomology classes 
$\gamma_{1}, \ldots, \gamma_{k} \in H^{*}(V, {\Bbb C})$ and non-negative integers
$n_{1}, \ldots, n_{k}$, the corresponding descendant Gromov-Witten invariants
are defined by
\[ \left< \tau_{n_{1}}(\gamma_{1}) \cdots
        \tau_{n_{k}}(\gamma_{k}) \right>_{g, A}  
	= \int_{\left[\overline{\cal M}_{g, k}(V, A) \right]^{\rm virt}}
        c_{1}(E_{1})^{n_{1}} \cup {\rm ev}_{1}^{*}(\gamma_{1}) \cup
                \cdots \cup
        c_{1}(E_{k})^{n_{k}} \cup {\rm ev}_{k}^{*}(\gamma_{k}),
\]
where $\left[\overline{\cal M}_{g, k}(V, A) \right]^{\rm virt}$
 is the virtual fundamental class of $\overline{\cal M}_{g, k}(V, A)$
(cf. \cite{LiT1}). When all $n_{i}$'s are zero, the corresponding
invariants are called primary Gromov-Witten invariants.
The notation $\tau_{n, \alpha}$ which was used in the introduction
will be explained in the next subsection.

\subsection{Convention of notations}
\label{sec:notation}

We will use $d$ to denote the complex dimension of $V$ and
let $N$ be the dimension of the space of cohomology classes
$H^{*}(V, {\Bbb C})$.
To define the generating functions, we need to fix a basis
$\{\gamma_{1}, \ldots, \gamma_{N}\}$ of $H^{*}(V, {\Bbb C})$ 
with $\gamma_{1}$ equal to the
identity of the cohomology ring of $V$ and 
$\gamma_{\alpha} \in H^{p_{\alpha}, q_{\alpha}}(V, {\Bbb C})$ for every
$\alpha$. 
We also arrange the basis in such a way  that the dimension
of $\gamma_{\alpha}$ is non-decreasing with respect to $\alpha$
and if two cohomology classes have the same dimension, we also require that 
the holomorphic dimension $p_{\alpha}$ is non-decreasing.
We will abbreviate $\tau_{n}(\gamma_{\alpha})$ as
$\tau_{n, \alpha}$ and identify $\tau_{0, \alpha}$ with $\gamma_{\alpha}$. 
For each $\tau_{n, \alpha}$, we associate a parameter
$t_{n}^{\alpha}$ and the collection of all such parameters is
denoted by 
$T = (t_{n}^{\alpha} \mid n \in {\Bbb Z}_{+}, \, \alpha= 1, \ldots, N)$,
where ${\Bbb Z}_{+}$ is the set of non-negative integers. 
The space of all $T$'s is the big phase space and its subspace
$\{T \mid t_{n}^{\alpha} = 0 \, \, {\rm if} \, \, n > 0 \}$ is the small
phase space.  For convenience, we will always identify the symbol 
$\tau_{n, \alpha}$ with the tangent vector field 
$\frac{\partial}{\partial t^{\alpha}_{n}}$ on the big phase space.
We also consider $\tau_{n, \alpha}$ with $n<0$ as a zero operator.
On the small phase space, we write $t_{0}^{\alpha}$ simply as
$t^{\alpha}$ and also identify the cohomology class $\gamma_{\alpha}$
with the vector field $\frac{\partial}{\partial t^{\alpha}}$.

As in the introduction, we can define the partition function $Z$ and
free energy function $F_{g}$ on the big phase
space. These are the generating functions of the corresponding classes
of Gromov-Witten invariants. The restriction of $F_{g}$ to the small phase 
space are denoted by $F_{g}^{s}$.
As in \cite{LT}, we will denote the tensor defined by the $k$-th
covariant derivative of $F_{g}$ by
$\left<\left< \right. \right. \underbrace{\cdot \cdots \cdot}_{k} 
 	\left. \left. \right>\right>_{g}$.
This is a symmetric $k$-tensors on the big phase
space defined by
\[ \left<\left< \tau_{m_{1}, \alpha_{1}}
        \tau_{m_{2}, \alpha_{2}}
        \cdots \tau_{m_{k}, \alpha_{k}}
         \right>\right>_{g} := 
        \frac{\partial^{k}}{\partial t^{\alpha_{1}}_{m_{1}}
        \partial t^{\alpha_{2}}_{m_{k}} \cdots
        \partial t^{\alpha_{k}}_{m_{k}}} F_{g}.
\]
This tensor is called the {\it $k$-point (correlation) function}.
The corresponding tensor on the small phase space is denoted by
$\left<\left< \right.\right.  \underbrace{\cdot \cdots \cdot}_{k} 
	\left. \left. \right>\right>_{g, s}$.

Besides the above notations, 
we will also use the following convention
throughout the paper unless otherwise stated.
Lower case Greek letters, e.g. 
$\alpha$, $\beta$, $\mu$, $\nu$, $\sigma$,..., etc., will
be used to index the cohomology classes. The range of these indices is from
$1$ to $N$, where $N$ is the dimension of the space of cohomology classes. 
Lower case English letters, e.g. $i$, $j$, $k$, $m$, $n$,
..., etc., will be used to index the level of descendents. Their range is the
set of all non-negative integers, i.e. ${\Bbb Z}_{+}$.
All summations are over the entire ranges of the indices unless otherwise
indicated. 
Let $\eta_{\alpha \beta} = \int_{V} \gamma_{\alpha} \cup \gamma_{\beta}$
be the intersection form on $H^{*}(V, {\Bbb C})$. 
We will use $\eta = (\eta_{\alpha \beta})$ and 
$\eta^{-1} = (\eta^{\alpha \beta})$ to lower and raise indices. Let
${\cal C} = ({\cal C}_{\alpha}^{\beta})$ be the matrix of multiplication
by the first Chern class $c_{1}(V)$ in the ordinary cohomology ring, i.e.
\begin{equation} \label{eqn:chernmatrix}
 c_{1}(V) \cup \gamma_{\alpha} = \sum_{\beta} {\cal C}_{\alpha}^{\beta}
                \gamma_{\beta}. 
\end{equation}
Since we are dealing with even dimensional cohomology classes only,
both $\eta$ and ${\cal C} \eta$ are symmetric matrices, where the entries
of ${\cal C}\eta$ are given by
${\cal C}_{\alpha \beta} = \int_{V}c_{1}(V) \cup \gamma_{\alpha} 
        \cup \gamma_{\beta}$.
Let 
\begin{equation} \label{eqn:balpha}
         b_{\alpha} = p_{\alpha} - \frac{1}{2}(d-1),
\end{equation}
where $d = {\rm dim}_{\Bbb C}V$.
The following simple observations will be used throughout the calculations
without mentioning: If $\eta^{\alpha \beta} \neq 0$ or
$\eta_{\alpha \beta} \neq 0$, then $b_{\alpha} = 1-b_{\beta}$.
${\cal C}_{\alpha}^{\beta} \neq 0$ implies $b_{\beta} =
1 + b_{\alpha}$, and ${\cal C}_{\alpha\beta} \neq 0$ implies 
$b_{\beta} = - b_{\alpha}$.

Instead of coordinates 
$\{t^{\alpha}_{m} \mid m\in {\Bbb Z}_{+}, \, \alpha = 1, \ldots, N\}$, 
it is very convenient to use the following shifted
coordinates on the big phase space
\begin{equation} \label{eqn:ttilde}
        \tilde{t}^{\alpha}_{m} = 
        t^{\alpha}_{m} - \delta_{m, 1} \delta_{\alpha, 1}
        = \left\{ \begin{array}{ll}
                        t^{\alpha}_{m}-1, & \textrm{ if }m=\alpha=1, \\
                        t^{\alpha}_{m}, & \textrm{ otherwise.}
                  \end{array}
                \right. 
\end{equation}

\subsection{Topological recursion relation}
\label{sec:TRR}

Topological recursion relations reduce the levels of descendants
in correlation functions.
The genus-$0$ {\it topological recursion relation} has the following
form (cf \cite{RT2} and \cite{W2}):
\[ \left<\left< \tau_{m, \alpha}
                \tau_{n, \beta}
                \tau_{k, \mu} \right>\right>_{0} =
        \sum_{\sigma} \left<\left< \tau_{m-1, \alpha}
                        \gamma_{\sigma} \right>\right>_{0}
                \left<\left< \gamma^{\sigma} 
                \tau_{n, \beta}
                \tau_{k, \mu} \right>\right>_{0}, \]
for $m > 0$. In this formula, we used the convention that the indices
of cohomology classes are raised by $\eta^{-1}$. Therefore $\gamma^{\sigma}$
should be understood as $\sum_{\rho}\eta^{\sigma \rho}\gamma_{\rho}$.
This recursion relation implies the following
{\it genus-0 constitutive relation} \cite{DW},
\[ \left<\left< \tau_{m, \alpha} \tau_{n, \beta} \right>\right>_{0}
	= \left< \tau_{m, \alpha} \tau_{n, \beta}  
		e^{\sum_{\sigma} u^{\sigma} \gamma_{\sigma}}\right>_{0}, \]
where $u^{\sigma} = \left<\left< \gamma_{1} \gamma^{\sigma} 
                        \right>\right>_{0}$.
This relation is an important building block in defining the $\tau$-function
for Frobenius manifolds \cite{Du}
(which corresponds to $F_{0}$ in the topological sigma model).

As noted by Witten \cite{W2}, the genus-0 topological recursion relation
implies the {\it generalized WDVV equation}:
\[ \sum_{\sigma} \left<\left< \tau_{m, \alpha}
                \tau_{n, \beta}
                \gamma_{\sigma} \right>\right>_{0}
                \left<\left< \gamma^{\sigma} 
                \tau_{k, \mu}
                \tau_{l, \nu} \right>\right>_{0}  
         = \sum_{\sigma} \left<\left< \tau_{m, \alpha}
                        \tau_{k, \mu} 
                        \gamma_{\sigma} \right>\right>_{0}
                \left<\left< \gamma^{\sigma} 
                \tau_{n, \beta}
                \tau_{l, \nu} \right>\right>_{0}. 
\]
When restricted to the small phase space, this equation is
usually called the {\it WDVV} equation.
It gives the associativity for the quantum cohomology
which is defined by  the third derivatives of $F_{0}^{s}$ and $\eta^{-1}$
(see Section~\ref{sec:Quantum}).
On the big phase space, this equation is the key ingredient 
in the proof of the genus-0 Virasoro conjecture (cf. \cite{LT}).

The {\it genus-1 topological recursion relation} is the following
\begin{equation} \label{eqn:g1TRR}
 \left<\left< \tau_{m+1, \alpha} \right>\right>_{1} =
        \sum_{\sigma} \left<\left< \tau_{m, \alpha}
                        \gamma_{\sigma} \right>\right>_{0}
                \left<\left< \gamma^{\sigma} \right>\right>_{1}
	+ \frac{1}{24} \sum_{\sigma} 
		\left<\left< \tau_{m, \alpha}
                    \gamma_{\sigma} \gamma^{\sigma} \right>\right>_{0}.
\end{equation}
This formula implies the {\it genus-1 constitutive relation} \cite{DW}
\begin{equation} \label{eqn:g1Cons}
 F_{1} = \left< e^{\sum_{\alpha} u^{\alpha} \gamma_{\alpha}} \right>_{1}
	+ \frac{1}{24} \log \det \left(\frac{\partial u^{\alpha}}
				{\partial t_{0}^{\beta}} \right),
\end{equation}
where $u^{\alpha} = \left<\left< \gamma_{1} \gamma^{\alpha} 
			\right>\right>_{0}$.

\subsection{Some special vector fields on the big phase space}
\label{sec:specialVF}
In \cite{LT}, we introduced several special vector fields
on the big phase space. These vector fields played very important
role in the proof of the genus-0 Virasoro conjecture.
The first one is the {\it string vector field}:
\[ {\cal S} := - \sum_{m, \alpha} \tilde{t}^{\alpha}_{m} 
         \tau_{m-1, \alpha}. \]
The restriction of ${\cal S}$ to the small phase space is just $\gamma_{1}$.
The famous {\it string equation} (cf. \cite{RT2} and \cite{W2})
can be expressed as
\[ \left<\left< {\cal S} \right>\right>_{g}  =
        \frac{1}{2} \delta_{g, 0} \sum_{\alpha, \beta} \eta_{\alpha \beta} 
                t^{\alpha}_{0} t^{\beta}_{0}. \]
This equation is equivalent to Eguchi, Hori, and Xiong's $L_{-1}$ 
constraint. 

The second vector field is the {\it Dilaton vector field}:
\[ {\cal D} := - \sum_{m, \alpha} \tilde{t}^{\alpha}_{m}
        \tau_{m, \alpha}. \]
When restricted to the small phase space, this vector field
does not tangent to the small phase space.
The so called {\it dilaton equation} is the following:
\[ \left<\left< {\cal D} \right>\right>_{g}  =
        (2g-2) F_{g} + \frac{1}{24} \, \chi(V)\delta_{g, 1}, \]
where $\chi(V)$ is the Euler characteristic number of $V$.
This equation implies the following (Lemma 1.2 in \cite{LT}):
\begin{equation} \label{eqn:Dilaton23}
\left<\left< {\cal D} \tau_{m, \alpha} \right>\right>_{0} =
                -\left<\left< \tau_{m, \alpha} \right>\right>_{0},
  \, \, \, \, \, {\rm and} \, \, \, \, \,
\left<\left< {\cal D} \tau_{m, \alpha} 
                \tau_{n, \beta} \right>\right>_{0} \equiv 0.
\end{equation}

Of particular importance is the following vector field:
\[ {\cal X} := - \sum_{m, \alpha} \left(m + b_{\alpha} - b_{1} - 1
                        \right)\tilde{t}^{\alpha}_{m} \tau_{m, \alpha} 
        - \sum_{m, \alpha, \beta} 
        {\cal C}_{\alpha}^{\beta}\tilde{t}^{\alpha}_{m}\tau_{m-1, \beta}.
        \]
When restricted to the small phase space, this vector field is the
Euler vector field $E$ mentioned in the introduction.
Therefore we also call ${\cal X}$ itself the {\it Euler vector field
(on the big phase space)}.
As noted in \cite{EHX1}, the divisor equation for
the first Chern class $c_{1}(V)$ together with the selection rule implies
the following {\it quasi-homogeneity equation}:
\[ \left<\left< {\cal X} \right>\right>_{g}  =
        2(b_{1}+1)(1-g) F_{g} 
        + \frac{1}{2} \delta_{g, 0}\sum_{\alpha, \beta} {\cal C}_{\alpha \beta}
        t^{\alpha}_{0}t^{\beta}_{0}
        - \frac{1}{24} \delta_{g, 1} \int_{V} c_{1}(V) \cup 
		c_{d-1}(V), \]
where $d$ is the complex dimension of $V$ and $c_{i}$ is the $i$-th
Chern class.
This equation implies the following (Lemma 1.4 in \cite{LT})
\begin{lem} \label{lem:EulerCorr}
\begin{eqnarray*}
&{\rm (i)}&  \left<\left< {\cal X} \right>\right>_{0} = 2(b_{1} + 1) F_{0}
                + \frac{1}{2} \sum_{\alpha, \beta} {\cal C}_{\alpha \beta}
        t^{\alpha}_{0}t^{\beta}_{0}.    \\
&{\rm (ii)}&   \left<\left< {\cal X} \tau_{m, \alpha} \right>\right>_{0} =
                \left(m + b_{\alpha} + b_{1} + 1 \right)
                \left<\left< \tau_{m, \alpha} \right>\right>_{0}
                + \sum_{\beta} {\cal C}_{\alpha}^{\beta}
        \left<\left< \tau_{m-1, \beta} \right>\right>_{0}  
    + \delta_{m, 0} \sum_{\beta} {\cal C}_{\alpha \beta} t^{\beta}_{0}.
                         \\
&{\rm (iii)} & \left<\left< {\cal X} \tau_{m, \alpha}
                \tau_{n, \beta} \right>\right>_{0} =
        \delta_{m, 0} \delta_{n,0} {\cal C}_{\alpha \beta} +
        (m+n+b_{\alpha}+b_{\beta}) \left<\left<\tau_{m, \alpha}
                \tau_{n, \beta} \right>\right>_{0} \\
 &  & \textrm{ \hspace{110pt}}  + \sum_{\mu} {\cal C}_{\alpha}^{\mu}
        \left<\left<\tau_{m-1, \mu}
                \tau_{n, \beta} \right>\right>_{0} 
        + \sum_{\mu} {\cal C}_{\beta}^{\mu}
        \left<\left<\tau_{m, \alpha}
                \tau_{n-1, \mu} \right>\right>_{0}.
\end{eqnarray*} 
\end{lem}

In \cite{LT}, we also introduced a sequence of vector fields
${\cal L}_{n}$ which are
the first derivative part of the Virasoro operators. The first four vector
fields are 
\begin{eqnarray} 
{\cal L}_{-1} & := & - {\cal S},    \nonumber \\
{\cal L}_{0} & := & - {\cal X} - (b_{1}+1){\cal D}, \nonumber \\
{\cal L}_{1} & := & \sum_{m, \alpha} (m+b_{\alpha})(m+b_{\alpha}+1) 
        \tilde{t}^{\alpha}_{m} 
                \tau_{m+1, \alpha}  
               \nonumber \\
        && + \sum_{m, \alpha , \beta} (2m+2b_{\alpha}+1) 
                {\cal C}_{\alpha}^{\beta}
        \tilde{t}^{\alpha}_{m} 
                \tau_{m, \beta} 
         + \sum_{m, \alpha, \beta} 
                ({\cal C}^{2})_{\alpha}^{\beta}
        \tilde{t}^{\alpha}_{m} \tau_{m-1, \beta}	\nonumber \\
{\cal L}_{2} & := & 
	 \sum_{m, \alpha} (m+b_{\alpha})(m+b_{\alpha}+1) 
                (m+b_{\alpha}+2) 
        \tilde{t}^{\alpha}_{m} \tau_{m+2, \alpha}  \nonumber \\
        && + \sum_{m, \alpha, \beta} 
        \left\{ 3(m+b_{\alpha})^{2} + 6(m+b_{\alpha}) +2 \right\} 
     {\cal C}_{\alpha}^{\beta}
        \tilde{t}^{\alpha}_{m} \tau_{m + 1, \beta} \nonumber \\
        && + \sum_{m, \alpha, \beta} 3(m+b_{\alpha}+1) 
                ({\cal C}^{2})_{\alpha}^{\beta}
        \tilde{t}^{\alpha}_{m} \tau_{m, \beta}  
         + \sum_{m, \alpha, \beta} 
                ({\cal C}^{3})_{\alpha}^{\beta}
        \tilde{t}^{\alpha}_{m} \tau_{m-1, \beta} \label{eqn:VirVectDef}
\end{eqnarray}
The following formulas were proved in \cite{LT} and will be used later:
\begin{eqnarray}
\left<\left< \gamma_{\mu}{\cal L}_{0} \gamma_{\nu} \right>\right>_{0} 
	& = & 
	- \left<\left< \gamma_{\mu}{\cal X} \gamma_{\nu} \right>\right>_{0} 
	 \nonumber \\
\left<\left< \gamma_{\mu}{\cal L}_{1} \gamma_{\nu} \right>\right>_{0}
	& = & - \sum_{\alpha}
	\left<\left< \gamma_{\mu}{\cal X} \gamma_{\alpha} \right>\right>_{0} 
	\left<\left< \gamma^{\alpha}{\cal X} \gamma_{\nu} \right>\right>_{0}
		\nonumber  \\
        && + \sum_{\alpha} b_{\alpha}(b_{\alpha}-1) 
            \left<\left< \gamma_{\alpha} \right>\right>_{0}
                \left<\left< \gamma^{\alpha}
                \gamma_{\mu} \gamma_{\nu})  \right>\right>_{0} 
		 \nonumber \\
\left<\left< \gamma_{\mu}{\cal L}_{2} \gamma_{\nu} \right>\right>_{0} 
	& = & - \sum_{\alpha, \beta}
	\left<\left< \gamma_{\mu}{\cal X} \gamma_{\alpha} \right>\right>_{0} 
	\left<\left< \gamma^{\alpha}{\cal X} \gamma_{\beta} \right>\right>_{0}
	\left<\left< \gamma^{\beta}{\cal X} \gamma_{\nu} \right>\right>_{0}
		\nonumber  \\
        && + \sum_{\alpha, \beta} b_{\alpha}(b_{\alpha}-1) 
            \left<\left< \gamma_{\alpha} \right>\right>_{0}
                \left<\left< \gamma^{\alpha}
                \gamma_{\mu} \gamma_{\beta})  \right>\right>_{0}
	\left<\left< \gamma^{\beta}{\cal X} \gamma_{\nu} \right>\right>_{0} 
		 \nonumber \\
        &&+ \sum_{\alpha} (b_{\alpha}-1) b_{\alpha}
                (b_{\alpha}+1) 
         \left<\left< \tau_{1, \alpha}
                \right>\right>_{0}  
                \left<\left< \gamma^{\alpha} 
                \gamma_{\mu} \gamma_{\nu} 
                \right>\right>_{0} \nonumber \\
        && + \sum_{\alpha, \beta} 
          (3 b_{\beta}^{2} - 1)
             {\cal C}_{\beta}^{\alpha}
        \left<\left< \gamma_{\alpha}
        \right>\right>_{0} 
                \left<\left< \gamma^{\beta} 
                \gamma_{\mu}
                \gamma_{\nu}) 
                \right>\right>_{0} \nonumber \\
        &&- \sum_{\alpha, \beta} (b_{\alpha}-1)b_{\alpha}(b_{\beta}-1)
         \left<\left< \gamma^{\alpha}
                \right>\right>_{0} 
         \left<\left< \gamma_{\alpha} \gamma_{\beta}
                \right>\right>_{0} 
             \left<\left< \gamma^{\beta} 
                \gamma_{\mu} \gamma_{\nu} 
                \right>\right>_{0} \nonumber \\
        && -    \sum_{\alpha, \beta}
                b_{\beta}(b_{\beta}+1){\cal C}_{\beta \alpha}
                \left<\left< \gamma^{\alpha}
                \right>\right>_{0} 
                \left<\left< \gamma^{\beta} 
                \gamma_{\mu} \gamma_{\nu}
                \right>\right>_{0}. \label{eqn:VirVect3pt}
\end{eqnarray} 
Due to Lemma 1.2 (3) in \cite{LT}, 
the first formula is just the definition of ${\cal L}_{0}$.
The second formula is a special case of the formula (19) 
in \cite{LT} plus the generalized WDVV equation. 
The third formula is a special case of the formula (26) 
in \cite{LT} plus the generalized WDVV equation and the second
formula.
Together with the obvious relation that the restriction of 
${\cal L}_{-1}$ to the small phase space is $- E^{0}$, these
formulas reveal an interesting relationship between the Virasoro operators
and the quantum powers of the Euler vector fields.
In fact, when restricted to the small phase space,
The first lines of the right hand sides of the above equations
are respectively 
$- \left<\left< \gamma_{\mu} E \gamma_{\nu} \right>\right>_{0, s}$,
$- \left<\left< \gamma_{\mu} E^{2} \gamma_{\nu} \right>\right>_{0, s}$,
$- \left<\left< \gamma_{\mu} E^{3} \gamma_{\nu} \right>\right>_{0, s}$.
With a slight modification of ${\cal L}_{n}$, the extra terms on
the right hand sides of the above equations may disappear.
This can be done by simply moving the extra terms to the left
hand sides, expressing them as 3-point functions with two arguments
equal to $\gamma_{\mu}$ and $\gamma_{\nu}$, then adding the third
arguments (which are again vector fields) to the corresponding
${\cal L}_{n}$'s (see also \cite{G2}).
We note here that for the second term on the right hand side of
the third equation, we can interchange the position of 
$\gamma_{\mu}$ and ${\cal X}$ (by the generalized WDVV equation), then using
Lemma~\ref{lem:EulerCorr} (iii) to remove ${\cal X}$. The third equation
 can then be simplified as
\begin{eqnarray}
\left<\left< \gamma_{\mu}{\cal L}_{2} \gamma_{\nu} \right>\right>_{0} 
	& = & - \sum_{\alpha, \beta}
	\left<\left< \gamma_{\mu}{\cal X} \gamma_{\alpha} \right>\right>_{0} 
	\left<\left< \gamma^{\alpha}{\cal X} \gamma_{\beta} \right>\right>_{0}
	\left<\left< \gamma^{\beta}{\cal X} \gamma_{\nu} \right>\right>_{0}
		\nonumber  \\
        &&+ \sum_{\alpha} (b_{\alpha}-1) b_{\alpha}
                (b_{\alpha}+1) 
         \left<\left< \tau_{1, \alpha}
                \right>\right>_{0}  
                \left<\left< \gamma^{\alpha} 
                \gamma_{\mu} \gamma_{\nu} 
                \right>\right>_{0} \nonumber \\
        &&+ \sum_{\alpha} (b_{\alpha}-1) b_{\alpha}
                (b_{\alpha}+1) 
         \left<\left< \gamma^{\alpha}
                \right>\right>_{0}  
                \left<\left< \tau_{1, \alpha} 
                \gamma_{\mu} \gamma_{\nu} 
                \right>\right>_{0} \nonumber \\
        && + \sum_{\alpha, \beta} 
          (3 b_{\beta}^{2} - 1)
             {\cal C}_{\beta}^{\alpha}
        \left<\left< \gamma_{\alpha}
        \right>\right>_{0} 
                \left<\left< \gamma^{\beta} 
                \gamma_{\mu}
                \gamma_{\nu}) 
                \right>\right>_{0}.  \label{eqn:3ptL2}
\end{eqnarray}

\subsection{Virasoro operators}
\label{sec:VirOp}
The first four Virasoro operators constructed by
Eguchi, Hori, and Xiong  are the following
\begin{eqnarray}
        L_{-1} &:= & {\cal L}_{-1}
                        + \frac{1}{2 \lambda^{2}} \sum_{\alpha, \beta}
                        \eta_{\alpha \beta} t^{\alpha}_{0} t^{\beta}_{0}, 
                        \nonumber \\
        L_{0} &:= & {\cal L}_{0}
		+ \frac{1}{2 \lambda^{2}} \sum_{\alpha, \beta}
                        {\cal C}_{\alpha \beta} t^{\alpha}_{0} t^{\beta}_{0}
               + \frac{1}{24}\left( (b_{1}+1) \chi(V) -
                        \int_{V} c_{1}(V) \cup c_{d-1}(V) \right),
		\nonumber \\
        L_{1} &:= & {\cal L}_{1}
		+ \frac{\lambda^{2}}{2}\sum_{\alpha} b_{\alpha}(1-b_{\alpha}) 
                 \gamma_{\alpha}\gamma^{\alpha}
		 + \frac{1}{2 \lambda^{2}}\sum_{\alpha, \beta} 
                ({\cal C}^{2})_{\alpha \beta} t^{\alpha}_{0} t^{\beta}_{0}
                 \nonumber \\
        L_{2} &:= & {\cal L}_{2}
		 - \lambda^{2}
		\sum_{\alpha} (b_{\alpha}-1) b_{\alpha}(b_{\alpha}+1) 
           	\tau_{1, \alpha} \gamma^{\alpha}
        	- \frac{\lambda^{2}}{2}\sum_{\alpha, \beta} 
		(3b_{\alpha}^{2} -1)
                {\cal C}_{\alpha}^{\beta} 
                 \gamma_{\beta}\gamma^{\alpha}
		 \nonumber \\
           	&& \hspace{50pt}+ \frac{1}{2 \lambda^{2}}\sum_{\alpha, \beta} 
                ({\cal C}^{3})_{\alpha \beta} t^{\alpha}_{0} t^{\beta}_{0}.
\end{eqnarray}
Because of the Virasoro relation
\[      [L_{m}, L_{n}] = (m-n)L_{m+n}, \]
for $m$, $n \geq -1$, the above operators generate all $L_{n}$ operators
with $n \geq -1$. We will not consider $L_{n}$ operators with $n < -1$
in this paper. Recall that the $L_{-1}$-constraint is equivalent to the string
equation, which is valid for all manifolds.
Due to the above Virasoro bracket relation, to prove the Virasoro conjecture,
it suffices to prove the $L_{2}$-constraint. 

\subsection{Review of the genus-0 Virasoro conjecture}
\label{sec:g0Vir}
In \cite{EHX2}, a heuristic argument for deriving the genus-$0$ 
constraints for ${\Bbb C}P^{n}$ was given.
It was pointed out in \cite{LT} that there is a serious gap in
this derivation. This observation was confirmed
by conversations between authors of \cite{LT}
and authors of \cite{EHX2} both before and after the paper 
\cite{LT} was written.

The first complete proof of the genus-0 Virasoro conjecture was
given in \cite{LT}. Actually the conjecture posed in \cite{EHX2}
was only for Fano varieties with vanishing Hodge numbers $h^{p,q}(V, {\Bbb C})$
for $p \neq q$ (cf. \cite{Bor}). This was later improved to cover all
compact smooth K\"{a}hler manifolds in \cite{EJX} after the Virasoro
operators were modified according to a suggestion of S. Katz.
It was pointed out for the first time in \cite{LT}
that the genus-0 Virasoro conjecture does
not need any assumption on manifolds, i.e. it is also valid for all compact
symplectic manifolds (The modification suggested by Katz does not apply to
general symplectic manifolds since it involves holomorphic dimensions of 
cohomology classes). Besides the
two constraints known before, i.e. $L_{-1}$ and $L_{0}$ constraints,
the key ingredients used in \cite{LT} were the genus-0
topological recursion relation and the generalized WDVV equation (those
equations are also valid for all Frobenius manifolds \cite{Du}). More
precisely, we computed the following expression
\[        \sum_{\sigma}
        \left<\left< {\cal L}_{n} \left({\cal L}_{0} - (n+1){\cal D}\right)
                \gamma_{\sigma} \right>\right>_{0} 
        \left<\left< \gamma^{\sigma} \tau_{k, \mu}
                \tau_{l, \nu} 
                \right>\right>_{0} 
        -  \left<\left< {\cal L}_{n} \tau_{k, \mu}
                \gamma_{\sigma} \right>\right>_{0}
        \left<\left< \gamma^{\sigma} \left({\cal L}_{0} - (n+1){\cal D}\right)
                \tau_{l, \nu}
                \right>\right>_{0}.  
\]
On the one hand, by the generalized WDVV equation, this expression
is 0. On the other hand, if the $L_{n}$ constraint is correct, we can
compute each 3-point function in the above expression separately, and
when combining the results together and using the 
genus-0 topological recursion relation, we can show that this expression
is just $\frac{\partial}{\partial t^{\nu}_{l}} 
        \frac{\partial}{\partial t^{\mu}_{k}} \Psi_{0, n+1}$. Once
we know that all second derivatives of $\Psi_{0, n+1}$ are zero,
the dilaton equation then trivially implies that $\Psi_{0, n+1} =0$.
This gives an inductive proof to the genus-0 Virasoro conjecture.
The key point in this proof is to observe the above relation between
the generalized WDVV equation and the Virasoro conjecture. Once this
relation is observed, the computation involved are quite straightforward,
although it is a little tedious. In Section 3 and 4 of \cite{LT},
we gave the full details of the computations. The advantage for doing so
instead of giving a more concise presentation is that one can
see clearly how each term of these complicated operators evolves during this
process. In particular, one can see how terms of $L_{n+1}$ emerge
from the above manipulation of expressions involving only ${\cal L}_{n}$
and ${\cal D}$.

In \cite{LT}, the same method  was also used to give
the first proof to another sequence of genus-0 constraints, called
$\widetilde{L}_{n}$ constraints, which were also conjectured in \cite{EHX2}.
Note that in the derivation of \cite{EHX2}, the two sequences of
constraints ($\{L_{n}\}$ and  $\{\widetilde{L}_{n}\}$ constraints)
are always mingled together. It is not clear how to separate these 
two sequences using the original arguments in \cite{EHX2}. 
To complete the proof along
the original lines of \cite{EHX2}, 
 one  needs to prove 
$\tilde{L}_{1}$ and $\tilde{L}_{2}$ constraints
first
(In the presentation of \cite{G2}, it is not
clear what kind of role the $\tilde{L}_{2}$ constraint plays, while in 
\cite{EHX2} this constraint was mixed with the $L_{2}$ constraint.).
It was noticed for the first time in \cite{LT} that these two sequences 
can be treated completely independently. The method for proving them
are the same. If one knows
how to prove one sequence, one also knows how to handle another one.

After \cite{LT} was submitted to journal and posted on the web, 
an alternative proof to the genus-0
Virasoro conjecture was given in \cite{DZ2}. In fact, 
the genus-0 Virasoro constraints were extended in \cite{DZ2} 
to the setting of abstract
Frobenius manifolds, which are defined by solutions of the WDVV equation
and also by axiomizing basic properties of the quantum cohomology. 
Since in genera bigger than 1, the corresponding
constitutive relations do not exist yet, it is not clear how to define
the analogue of $F_{g}$ for abstract Frobenius manifolds. Therefore it is
not clear how to interpret the full Virasoro conjecture for this setting.
The third proof to the genus-0 Virasoro conjecture
was given in \cite{G2} by combining arguments in \cite{EHX2} and \cite{DZ2}.

\subsection{Quantum cohomology}
\label{sec:Quantum}

At each point of the small phase space, which is identified with
$H^{*}(V, {\Bbb C})$, we can define a new product structure among
 cohomology classes, called the {\it Quantum product}, in the following way:
\[ \gamma_{\alpha} \bullet \gamma_{\beta} = \sum_{\sigma}
	\left<\left< \gamma_{\alpha} \gamma_{\beta} \gamma^{\sigma}
		\right>\right>_{0, s} \gamma_{\sigma}. \]
This product is commutative and associative (due to the WDVV equation).
In this way, we obtain new ring structures on $H^{*}(V, {\Bbb C})$,
which are called {\it quantum cohomologies of V}.
Since the restriction of the string vector field to the small phase
space is $\gamma_{1}$, the string equation implies the following
\begin{lem} \label{lem:StringCorr}
\[ \left<\left< \gamma_{1} \gamma_{\alpha}
                \gamma_{\beta} \right>\right>_{0, s} =
                 \eta_{\alpha \beta}, \,\,\,\,
{\rm and}\,\,\,\, \left<\left< \gamma_{1} \gamma_{\mu_{1}} \cdots
                \gamma_{\mu_{k}} \right>\right>_{0, s} = 0 \,\,\,\,
                {\rm if } \,\,\,\, k \geq 3. 
\]
\end{lem}
Especially the first equation in the lemma
tells us that $\gamma_{1}$ is always the identity of the quantum
cohomology no matter which point in the small phase space is chosen.

Since $H^{*}(V, {\Bbb C})$ is a linear space, we can identify tangent
spaces of $H^{*}(V, {\Bbb C})$ with $H^{*}(V, {\Bbb C})$ itself.
Therefore we can take quantum product for any two vector
fields on $H^{*}(V, {\Bbb C})$. 
The intersection form $\eta$ defines a flat metric (non-Riemannian)
on $H^{*}(V, {\Bbb C})$. Let $\nabla$ be the corresponding Levi-Civita
connection. It is straightforward to verify the following
\begin{equation} \label{eqn:derCorr}
 u \left<\left< v_{1} \cdots v_{k}
				\right>\right>_{g, s}
	=  \left<\left< u v_{1} \cdots v_{k}
				\right>\right>_{g, s}
	  + \sum_{i=1}^{k} \left<\left< v_{1} \cdots
		(\nabla_{u} v_{i}) \cdots v_{k}
				\right>\right>_{g, s},
\end{equation}
for any vector fields $u$, and $v_{1}, \ldots, v_{k}$ on the small
phase space. A simple application of this formula is the following
\begin{equation} \label{eqn:derProd}
	\nabla_{u}(v \bullet w) = (\nabla_{u}v) \bullet w 
				 + v \bullet (\nabla_{u} w)
			+ \sum_{\alpha} \left<\left< uvw \gamma^{\alpha}
				\right>\right>_{0, s} \gamma_{\alpha},
\end{equation}
for any vector fields $u$, $v$, and $w$
on the small phase space.

The most important vector field on the small phase space
 is the Euler vector field $E$ defined
in the introduction. It is the restriction to the small phase
space of the vector field ${\cal X}$ defined in Section~\ref{sec:specialVF}.
Therefore Lemma~\ref{lem:EulerCorr} implies the following
\begin{lem} \label{lem:sEulerCorr}
\begin{eqnarray*}
&{\rm (i)}&  \left<\left< E \right>\right>_{0, s} = 2(b_{1} + 1) F_{0}^{s}
                + \frac{1}{2} \sum_{\alpha, \beta} {\cal C}_{\alpha \beta}
        t^{\alpha}t^{\beta}.    \\
&{\rm (ii)}&   \left<\left< E \gamma_{\alpha} \right>\right>_{0, s} =
                \left(b_{\alpha} + b_{1} + 1 \right)
                \left<\left< \gamma_{\alpha} \right>\right>_{0, s}          
	    +  \sum_{\beta} {\cal C}_{\alpha \beta} t^{\beta}.
                        \\
&{\rm (iii)}& \left<\left< E \gamma_{\alpha}
                \gamma_{\beta} \right>\right>_{0, s} =
         {\cal C}_{\alpha \beta} +
        (b_{\alpha}+b_{\beta}) \left<\left<\gamma_{\alpha}
                \gamma_{\beta} \right>\right>_{0, s}. 
	\\
&{\rm (iv)}& \left<\left< E \gamma_{\alpha}
                \gamma_{\beta} \gamma_{\mu} \right>\right>_{0, s} =
         (b_{\alpha}+b_{\beta} + b_{\mu} - b_{1} -1) 
			\left<\left<\gamma_{\alpha}
                \gamma_{\beta} \gamma_{\mu} \right>\right>_{0, s}. 
	\\
&{\rm (v)}& \left<\left< E \gamma_{\alpha}
                \gamma_{\beta} \gamma_{\mu} \gamma_{\nu}\right>\right>_{0, s}
	 = (b_{\alpha}+b_{\beta} + b_{\mu} + b_{\nu} - 2b_{1} -2) 
			\left<\left<\gamma_{\alpha}
                \gamma_{\beta} \gamma_{\mu} \gamma_{\nu} 
		\right>\right>_{0, s}.
\end{eqnarray*} 
\end{lem}
A simple application of the third formula in this lemma is the following:
\begin{equation} \label{eqn:derEuler}
 v \left<\left< E \gamma_{\alpha}
                \gamma_{\beta} \right>\right>_{0, s} =
               (b_{\alpha}+b_{\beta}) \left<\left< v \gamma_{\alpha}
                \gamma_{\beta} \right>\right>_{0, s},
\end{equation}
where $v$ is any vector field on the small phase space.
Let $E^{i}$ be the $i$-th quantum power of $E$. Then 
\[  \left<\left<  \gamma^{\alpha} E^{i}
                \gamma_{\beta} \right>\right>_{0, s} =
	\sum_{\mu_{1}, \ldots,  \mu_{i-1}} \prod_{j=1}^{i}
		\left<\left<  \gamma^{\mu_{j-1}} E
                \gamma_{\mu_{j}} \right>\right>_{0, s},
\]
where $\mu_{0} = \alpha$ and $\mu_{i} = \beta$. 
Applying equation (\ref{eqn:derEuler}) to each factor, we obtain
\begin{eqnarray} \label{eqn:derPowerEuler}
  && E^{k} \left<\left<  \gamma^{\alpha} E^{i}
                \gamma_{\beta} \right>\right>_{0, s} \nonumber \\ 
	& = &
               (b_{\beta} - b_{\alpha} + i) 
		\left<\left<  \gamma^{\alpha} E^{k+i-1}
                \gamma_{\beta} \right>\right>_{0, s}
			 \nonumber \\ &&
		- \sum_{j=1}^{\min\{i, k\} - 1} \sum_{\mu} b_{\mu}
		\left<\left<  \gamma^{\alpha} E^{j}
                \gamma_{\mu} \right>\right>_{0, s}
		\left<\left<  \gamma^{\mu} E^{k+i-1-j}
                \gamma_{\beta} \right>\right>_{0, s}
			\nonumber \\ &&
		+ \sum_{j=1}^{\min\{i, k\} - 1} \sum_{\mu} b_{\mu}
		\left<\left<  \gamma^{\alpha} E^{k+i-1-j}
                \gamma_{\mu} \right>\right>_{0, s}
		\left<\left<  \gamma^{\mu} E^{j}
                \gamma_{\beta} \right>\right>_{0, s},
\end{eqnarray}
for $i$, $k \geq 1$.
For convenience, we will write 
\begin{equation} \label{eqn:EularCoef}
	 E^{k} = \sum_{\alpha} x_{k}^{\alpha} \gamma_{\alpha}, \,\,\,\,
	{\rm where} \, \, \, \, 
	 x_{k}^{\alpha} = \left<\left<  \gamma_{1} E^{k}
                \gamma^{\alpha} \right>\right>_{0, s}.
\end{equation}
Since $[E^{k}, \, E^{m}] = \sum_{\alpha} \left(E^{k} x_{m}^{\alpha}
			- E^{m} x_{k}^{\alpha}\right) \gamma_{\alpha}$,
a simple application of  equation (\ref{eqn:derPowerEuler}) proves
\[ [E^{k}, \, E^{m}] = (m-k) E^{m+k-1} \]
 for $m$, $k \geq 1$. If one of
$m$ and $k$ is equal to 0, the corresponding formula follows from
equation (\ref{eqn:derProd}) and Lemma~\ref{lem:StringCorr} since
$[ \gamma_{1}, \, E^{k}] = \nabla_{\gamma_{1}} E^{k}$. This gives
a simple proof to equation (\ref{eqn:VirEuler}).

\section{Relations between genus-0 and genus-1 data}
\label{sec:g0g1}

In this section, we will study  how much genus-1 information
can be obtained from genus-0 data. In particular, we will prove
Theorem~\ref{thm:EulerF1}. We first define two symmetric 4-tensors
$G_{0}$ and $G_{1}$ on the small phase space.
Let $S_{4}$ be the permutation group of 4 elements which acts
on the set $\{1, 2, 3, 4\}$. For any vector
fields $v_{1}, \ldots v_{4}$ on the small phase space, we define
\begin{eqnarray*}
G_{0}(v_{1}, v_{2}, v_{3}, v_{4}) & = &
	\sum_{g \in S_{4}} \sum_{\alpha, \beta} \left\{
		\frac{1}{6} \left<\left< v_{g(1)} v_{g(2)}v_{g(3)} 
					\gamma^{\alpha} \right>\right>_{0, s}
		\left<\left< \gamma_{\alpha} v_{g(4)} \gamma_{\beta}
					\gamma^{\beta} \right>\right>_{0, s}
		\right.	\\
	&& \hspace{40pt} 
	+ \frac{1}{24} \left<\left< v_{g(1)} v_{g(2)}v_{g(3)} v_{g(4)}
					\gamma^{\alpha} \right>\right>_{0, s}
	\left<\left< \gamma_{\alpha}  \gamma_{\beta}
					\gamma^{\beta} \right>\right>_{0, s}
			\\
	&& \hspace{40pt} \left.
	- \frac{1}{4} \left<\left< v_{g(1)} v_{g(2)}
			\gamma^{\alpha} \gamma^{\beta} \right>\right>_{0, s}
		\left<\left< \gamma_{\alpha}  \gamma_{\beta}
				v_{g(3)} v_{g(4)} \right>\right>_{0, s} 
			\right\},
\end{eqnarray*}
and 
\begin{eqnarray*}
G_{1}(v_{1}, v_{2}, v_{3}, v_{4}) & = &
	\sum_{g \in S_{4}} 
		3 \left<\left< \{v_{g(1)} \bullet v_{g(2)} \}
			\{v_{g(3)} \bullet v_{g(4)} \}
					 \right>\right>_{1, s}
		\\
	&&
	- \sum_{g \in S_{4}} 
	4 \left<\left< \{v_{g(1)} \bullet v_{g(2)} \bullet
			v_{g(3)} \} v_{g(4)} 
					 \right>\right>_{1, s} 
	\\
	&& 
	 - \sum_{g \in S_{4}} \sum_{\alpha}
		\left<\left< \{ v_{g(1)} \bullet v_{g(2)} \}
				v_{g(3)} v_{g(4)}
					\gamma^{\alpha} \right>\right>_{0, s}
	\left<\left< \gamma_{\alpha} \right>\right>_{1, s}
		\\
	&& 
	 + \sum_{g \in S_{4}} \sum_{\alpha}
		 2 \left<\left< v_{g(1)}  v_{g(2)} v_{g(3)} 
					\gamma^{\alpha} \right>\right>_{0, s}
	\left<\left< \{\gamma_{\alpha} \bullet v_{g(4)} \} 
			\right>\right>_{1, s}.	
\end{eqnarray*}
Note that $G_{0}$ is determined solely by genus-0 data, while each term in
$G_{1}$ contains genus-1 information. These two tensors are connected
by the following equation:
\begin{equation} \label{eqn:G}
	G_{0} + G_{1} = 0.
\end{equation}
This equation was proved in \cite{G1} where it was written in a different form.
The above formulation is a slight modification of the one given in \cite{DZ1}.
We first study the function $G_{1}$.
\begin{pro} \label{prop:G1}
\begin{eqnarray*}
	G_{1}(v_{1}, v_{2}, v_{3}, v_{4}) 
	&=& \sum_{g \in S_{4}} \left\{
		3 \{v_{g(1)} \bullet v_{g(2)} \} \left<\left< 
			\{v_{g(3)} \bullet v_{g(4)} \}
					 \right>\right>_{1, s}
		\right. \\
	&& \hspace{30pt}
	- 4 v_{g(4)} \left<\left< \{v_{g(1)} \bullet v_{g(2)} \bullet
			v_{g(3)} \} 
					 \right>\right>_{1, s} 
	\\
	&& \hspace{30pt} \left.
	 - 6 \left<\left< \left\{[ \{v_{g(1)} \bullet v_{g(2)} \}, \, 
			v_{g(3)}] \bullet v_{g(4)} \right\}
					 \right>\right>_{1, s} \right\}.
\end{eqnarray*}
\end{pro}
{\bf Proof:} Using equation (\ref{eqn:derCorr}) and (\ref{eqn:derProd})
to compute $\sum_{g \in S_{4}} 
		 3\{v_{g(1)} \bullet v_{g(2)} \} \left<\left< 
			\{v_{g(3)} \bullet v_{g(4)} \}
					 \right>\right>_{1, s}$
and $-\sum_{g \in S_{4}} 
  	4 v_{g(4)} \left<\left< \{v_{g(1)} \bullet v_{g(2)} \bullet
                        v_{g(3)} \}  \right>\right>_{1, s}$, then combining
the results together and using the symmetry of the tensors, we obtain the
desired formula. $\Box$

Applying this proposition to quantum powers of the Euler vector field $E$
and using equation (\ref{eqn:VirEuler}), we obtain the following
\begin{cor} \label{cor:G1Euler}
\begin{eqnarray*}
	&&	\frac{1}{24}G_{1}(E^{m_{1}}, E^{m_{2}}, E^{m_{3}}, E^{m_{4}})
		\\
	& = &  (2m_{1}+m) \left<\left< E^{m-1} \right>\right>_{1, s}
	+\sum_{i=2}^{4} 
		 E^{m_{1} + m_{i}} \left<\left< 
			E^{m-m_{1}- m_{i}}
					 \right>\right>_{1, s}
	- \sum_{i=1}^{4}  E^{m_{i}}\left<\left< E^{m-m_{i}} 
					 \right>\right>_{1, s},
\end{eqnarray*}
where $m_{1}, \ldots, m_{4}$ are arbitrary non-negative integers
and $m = m_{1}+m_{2}+m_{3}+m_{4}$.
\end{cor}
In the rest of this paper, we will use the following simple formulas
without mentioning:
\begin{lem} \label{lem:E0E1F1}
\[ \begin{array}{rl}
	{\rm (i)} & \left<\left< E^{0} \right>\right>_{1, s} = 0, \\
	{\rm (ii)} & \left<\left< E \right>\right>_{1, s} = - \frac{1}{24} 
		\int_{V} c_{1}(V) \cup c_{d-1}(V),  \\
	{\rm (iii)} & E^{0} \left<\left< E^{m} \right>\right>_{1, s} = 
		m \left<\left< E^{m-1} \right>\right>_{1, s}, \\
	{\rm (iv)} & E \left<\left< E^{m} \right>\right>_{1, s} = 
		(m-1) \left<\left< E^{m} \right>\right>_{1, s},
	\end{array} \]
for any non-negative integer $m$.
\end{lem}
{\bf Proof:} The first two equations are the restrictions of the genus-1
string equation and quasi-homogeneity equation to the small phase space
respectively. The last two equations follows from the first two equations
and equation (\ref{eqn:VirEuler}). $\Box$

A special case of the Corollary~\ref{cor:G1Euler} is the following
\begin{eqnarray} \label{eqn:G1EE}
 \frac{1}{24}G_{1}(E^{m-2-i}, E^{i}, E, E) 
	&=& E^{m-2} \left<\left< E^{2} \right>\right>_{1, s}
		-E^{m-i-2} \left<\left< E^{i+2} \right>\right>_{1, s}
		\nonumber \\ 
	&&    + 2 E^{m-i-1} \left<\left< E^{i+1} \right>\right>_{1, s}	
	   - E^{m-i} \left<\left< E^{i} \right>\right>_{1, s}  
\end{eqnarray}
for $1 \leq i \leq \left[\frac{m}{2}\right] -2 $ where 
$\left[\frac{m}{2}\right]$ is the largest integer which is less than
or equal to $\frac{m}{2}$.
If $m$ is even, Corollary~\ref{cor:G1Euler} implies
\begin{eqnarray*}
\frac{1}{48}G_{1}(E^{m/2-1}, E^{m/2-1}, E, E) & = &
	 \frac{1}{2} E^{m-2} \left<\left< E^{2} \right>\right>_{1, s}
	 - E^{m/2 + 1} \left<\left< E^{m/2-1} \right>\right>_{1, s}
		\nonumber \\ 
	&&   + E^{m/2} \left<\left< E^{m/2} \right>\right>_{1, s}
		- \left<\left< E^{m-1} \right>\right>_{1, s}	   
\end{eqnarray*}
Summing up equation (\ref{eqn:G1EE}) over 
$1 \leq i \leq \left[\frac{m}{2}\right] -2 $ and adding the above equation,
we obtain
\begin{eqnarray} \label{eqn:G1Even}
 && 	\frac{1}{48}G_{1}(E^{m/2-1}, E^{m/2-1}, E, E) 
	+ \sum_{i=1}^{m/2 -2}
	           \frac{1}{24}G_{1}(E^{m-2-i}, E^{i}, E, E) 
		\nonumber \\ 
	&=&  \frac{m-1}{2} E^{m-2} \left<\left< E^{2} \right>\right>_{1, s}  
	   - \left<\left< E^{m-1} \right>\right>_{1, s}
\end{eqnarray}
when $m$ is an even integer.
If $m$ is odd, Corollary~\ref{cor:G1Euler} implies
\begin{eqnarray*}
\frac{1}{24}G_{1}(E^{(m-1)/2}, E^{(m-3)/2}, E, E) & = &
	 E^{m-2} \left<\left< E^{2} \right>\right>_{1, s}
	 -  E^{(m+3)/2} \left<\left< E^{(m-3)/2} \right>\right>_{1, s}
		\nonumber \\ 
	&&   + E^{(m+1)/2} \left<\left< E^{(m-1)/2} \right>\right>_{1, s}
		- \left<\left< E^{m-1} \right>\right>_{1, s}	   
\end{eqnarray*}
Summing up equation (\ref{eqn:G1EE}) over 
$1 \leq i \leq \left[\frac{m}{2}\right] -2 $ and adding the above equation,
we obtain
\begin{equation} \label{eqn:G1Odd}
 \sum_{i=1}^{(m-3)/2}
	           \frac{1}{24}G_{1}(E^{m-2-i}, E^{i}, E, E) 
 =  \frac{m-1}{2} E^{m-2} \left<\left< E^{2} \right>\right>_{1, s}  
	   - \left<\left< E^{m-1} \right>\right>_{1, s}
\end{equation}
when $m$ is an odd integer.
Using the symmetry of the tensor $G_{1}$, 
we can express equation (\ref{eqn:G1Even}) and (\ref{eqn:G1Odd}) in a unified
form, which together with equation (\ref{eqn:G}) implies the following 
\begin{thm} \label{thm:g0g1G0} 
For an arbitrary manifold $V$,
\[ 
  \frac{m-1}{2} E^{m-2} \left<\left< E^{2} \right>\right>_{1, s}  
	   - \left<\left< E^{m-1} \right>\right>_{1, s}
	= - \sum_{i=1}^{m-3}
	           \frac{1}{48}G_{0}(E^{m-2-i}, E^{i}, E, E)
\]
for any integer $m \geq 2$.
\end{thm}
Since $G_{0}$ is defined by derivatives of $F_{0}^{s}$, 
this theorem in particular implies Theorem~\ref{thm:EulerF1}.

\section{A sequence of genus-0 functions}
\label{sec:phi-k}

Theorem~\ref{thm:g0g1G0} tells us that for $k \geq 3$,
$\left<\left< E^{k} \right>\right>_{1, s}$ 
can be computed in terms of
$\left<\left< E^{2} \right>\right>_{1, s}$ and some genus-0 data.
We will see later that the restriction of the genus-1 $L_{1}$ constraint
to the small phase space is equivalent to 
$ \left<\left< E^{2} \right>\right>_{1, s} = \phi_{2}$ where $\phi_{2}$
is defined in (\ref{eqn:phi2}). 
We can rewrite $\phi_{2}$ in the following form:
\begin{equation} \label{eqn:phi2New}
 \phi_{2} = - \frac{1}{24} \sum_{\alpha} 
	 \left<\left< E E \gamma_{\alpha} \gamma^{\alpha} \right>\right>_{0, s}
	+ \frac{1}{2}\sum_{\alpha}  \left(b_{\alpha} (1-b_{\alpha})
		 - \frac{b_{1}+1}{6}\right)
	\left<\left< \gamma_{\alpha} \gamma^{\alpha} \right>\right>_{0, s}.
\end{equation}
Motivated by Theorem~\ref{thm:g0g1G0}, we define
\begin{equation} \label{eqn:phi-k}
 \phi_{k} := \frac{k}{2} E^{k-1} \phi_{2} + 
	\sum_{i=1}^{k-2}
                   \frac{1}{48}G_{0}(E^{k-1-i}, E^{i}, E, E),
\end{equation}
for $k \geq 3$. For convenience, we also define
\begin{equation} \label{eqn:phi0phi1}
\phi_{0} := 0, \, \, \, \, {\rm and} \, \, \, \,
\phi_{1} := - \frac{1}{24} \int_{V} c_{1}(V) \cup c_{d-1}(V).
\end{equation}
The string equation and the quasi-homogeneity equation implies
\[
\phi_{0} = \left<\left< E^{0} \right>\right>_{1, s}, 
\, \, \, \, {\rm and} \, \, \, \,
\phi_{1} = \left<\left< E \right>\right>_{1, s}.
\]
An immediate consequence of Theorem~\ref{thm:g0g1G0} is the following
\begin{thm} \label{thm:EulerPhi}
For any manifold $V$, if 
$\left<\left< E^{2} \right>\right>_{1, s} = \phi_{2}$, then
$\left<\left< E^{k} \right>\right>_{1, s} = \phi_{k}$ for every $k$.
\end{thm}

The definition of $\phi_{k}$ given by (\ref{eqn:phi-k}) is hard to
use. For the convenience of later applications, we will give another
equivalent formulation in Theorem~\ref{thm:phik}. Before proving 
Theorem~\ref{thm:phik}, we need some preparations.
First, taking derivatives of the WDVV equation twice and
three times, we obtain the following
\begin{lem} \label{lem:derWDVV}
For any vector fields $u$, $v$, $w_{i}$ on the small phase space, we have
\begin{eqnarray*}
&{\rm (i)}& \sum_{\alpha}
	\left<\left< u w_{1} w_{2} \gamma_{\alpha} \right>\right>_{0, s}
	\left<\left< \gamma^{\alpha} v w_{3} \right>\right>_{0, s} \\
&& \hspace{20pt}
	+ \sum_{\alpha}
	\left<\left< u w_{1} \gamma_{\alpha} \right>\right>_{0, s}
	\left<\left< \gamma^{\alpha} v w_{2} w_{3} \right>\right>_{0, s} \\
&& =  \sum_{\alpha}
	\left<\left< v w_{1} w_{2} \gamma_{\alpha} \right>\right>_{0, s}
	\left<\left< \gamma^{\alpha} u w_{3} \right>\right>_{0, s} \\
&& \hspace{20pt}
	+ \sum_{\alpha}
	\left<\left< v w_{1} \gamma_{\alpha} \right>\right>_{0, s}
	\left<\left< \gamma^{\alpha} u w_{2} w_{3} \right>\right>_{0, s}, \\
&{\rm (ii)}& \sum_{\alpha}
	\left<\left< u w_{1} w_{2} \gamma_{\alpha} \right>\right>_{0, s}
	\left<\left< \gamma^{\alpha} v w_{3} w_{4} \right>\right>_{0, s} \\
&& \hspace{20pt}
	+ \sum_{\alpha}
	\left<\left< u w_{1} w_{3} \gamma_{\alpha} \right>\right>_{0, s}
	\left<\left< \gamma^{\alpha} v w_{2} w_{4} \right>\right>_{0, s} \\
&& \hspace{20pt}
	+ \sum_{\alpha}
	\left<\left< u w_{1} w_{2} w_{3} \gamma_{\alpha} \right>\right>_{0, s}
	\left<\left< \gamma^{\alpha} v w_{4} \right>\right>_{0, s} \\
&& \hspace{20pt}
	+ \sum_{\alpha}
	\left<\left< u w_{1}  \gamma_{\alpha} \right>\right>_{0, s}
	\left<\left< \gamma^{\alpha} v w_{2} w_{3} 
		w_{4} \right>\right>_{0, s} \\
&& = \sum_{\alpha}
	\left<\left< v w_{1} w_{2} \gamma_{\alpha} \right>\right>_{0, s}
	\left<\left< \gamma^{\alpha} u w_{3} w_{4} \right>\right>_{0, s} \\
&& \hspace{20pt}
	+ \sum_{\alpha}
	\left<\left< v w_{1} w_{3} \gamma_{\alpha} \right>\right>_{0, s}
	\left<\left< \gamma^{\alpha} u w_{2} w_{4} \right>\right>_{0, s} \\
&& \hspace{20pt}
	+ \sum_{\alpha}
	\left<\left< v w_{1} w_{2} w_{3} \gamma_{\alpha} \right>\right>_{0, s}
	\left<\left< \gamma^{\alpha} u w_{4} \right>\right>_{0, s} \\
&& \hspace{20pt}
	+ \sum_{\alpha}
	\left<\left< v w_{1}  \gamma_{\alpha} \right>\right>_{0, s}
	\left<\left< \gamma^{\alpha} u w_{2} w_{3} 
		w_{4} \right>\right>_{0, s}.
\end{eqnarray*}
\end{lem}
We can use these formulas and the WDVV equation to exchange positions
of two vector fields in a product of two correlation functions. Using
this lemma, we can prove the following
\begin{lem} \label{lem:G0munu}
For any $\mu$ and $\nu$,
\begin{eqnarray*}
G_{0}(\gamma_{\mu}, \gamma_{\nu}, E, E)
 &=& \sum_{\beta} \left<\left< E E (\gamma_{\mu} \bullet \gamma_{\nu}) 
		\gamma_{\beta} \gamma^{\beta} \right>\right>_{0, s} \\
 && + \sum_{\alpha, \beta} (2b_{\beta}-b_{\alpha}+b_{\mu} -1) 
	\left<\left< \gamma_{\mu} E \gamma^{\alpha} \right>\right>_{0, s}
	\left<\left< \gamma_{\alpha} \gamma_{\beta}
		\gamma^{\beta} \gamma_{\nu} \right>\right>_{0, s} \\
 && + \sum_{\alpha, \beta} (2b_{\beta}-b_{\alpha}+b_{\nu} -1) 
	\left<\left< \gamma_{\nu} E \gamma^{\alpha} \right>\right>_{0, s}
	\left<\left< \gamma_{\alpha} \gamma_{\beta}
		\gamma^{\beta} \gamma_{\mu} \right>\right>_{0, s} \\
 && + \sum_{\alpha, \beta} (b_{\alpha} - b_{1})
			(-4b_{\beta}-2b_{\alpha}+2b_{1} + 4) 
	\left<\left< \gamma_{\mu} \gamma_{\nu} \gamma^{\alpha} 
				\right>\right>_{0, s}
	\left<\left< \gamma_{\alpha} \gamma_{\beta}
		\gamma^{\beta} \right>\right>_{0, s} \\
 && + \sum_{\alpha, \beta} \left\{
	(-4b_{\beta}-2b_{\alpha}+b_{\mu} - b_{\nu} + 2)
	(b_{\mu} - b_{\nu} - 2 b_{\alpha} - 2 b_{\beta} + 2) \right. \\
 && \hspace{40pt}	
	\left. - 2 (b_{\nu} + b_{\alpha} + b_{\beta} - b_{1} - 1) \right\} 
	\left<\left< \gamma_{\mu}  \gamma^{\alpha} \gamma^{\beta}
			\right>\right>_{0, s}
	\left<\left< \gamma_{\alpha} \gamma_{\beta}
		 \gamma_{\nu} \right>\right>_{0, s}. 
\end{eqnarray*}
\end{lem}
{\bf Proof}: Applying Lemma~\ref{lem:derWDVV} (ii) with $u = E$, 
$v = \gamma^{\beta}$, $w_{1} = E$, $w_{2} = \gamma_{\mu}$, 
$w_{3} = \gamma_{\nu}$ and $w_{4} = \gamma_{\beta}$
to the expression
\begin{eqnarray*}
&& \sum_{\alpha, \beta}
        \left<\left< \gamma_{\mu} E E \gamma^{\alpha} \right>\right>_{0, s}
        \left<\left< \gamma_{\alpha} \gamma_{\nu}
		\gamma_{\beta} \gamma^{\beta} \right>\right>_{0, s} 
	+ \left<\left< \gamma_{\nu} E E \gamma^{\alpha} \right>\right>_{0, s}
        \left<\left< \gamma_{\alpha} \gamma_{\mu}
		\gamma_{\beta} \gamma^{\beta} \right>\right>_{0, s} \\
&& \hspace{50pt}
	+ \left<\left< \gamma_{\mu} \gamma_{\nu} 
			E E \gamma^{\alpha} \right>\right>_{0, s}
        \left<\left< \gamma_{\alpha}
		\gamma_{\beta} \gamma^{\beta} \right>\right>_{0, s}, 
\end{eqnarray*}
then applying Lemma~\ref{lem:derWDVV} (ii) again with $u = E$, 
$v = \gamma_{\mu}$, $w_{1} = E$, 
$w_{2} = \gamma_{\beta}$, 
$w_{3} = \gamma^{\beta}$ and $w_{4} = \gamma_{\nu}$
to the expression
\[ -\left\{\sum_{\alpha, \beta}
	\left<\left< E E \gamma^{\alpha} \right>\right>_{0, s}
        \left<\left< \gamma_{\alpha} \gamma_{\mu} \gamma_{\nu}
		\gamma_{\beta} \gamma^{\beta} \right>\right>_{0, s}
   	+ 2\left<\left< E E \gamma^{\beta}\gamma^{\alpha} \right>\right>_{0, s}
        \left<\left< \gamma_{\alpha} \gamma_{\beta} \gamma_{\mu} \gamma_{\nu}
		 \right>\right>_{0, s} \right\}, \]
after plugging the corresponding results into the definition of 
$G_{0}(\gamma_{\mu}, \gamma_{\nu}, E, E)$ and using
Lemma~\ref{lem:sEulerCorr} to 4-point and 5-point functions which
involve only one $E$, we obtain
\begin{eqnarray*}
G_{0}(\gamma_{\mu}, \gamma_{\nu}, E, E)
 &=& \sum_{\beta} \left<\left< E E (\gamma_{\mu} \bullet \gamma_{\nu}) 
                \gamma_{\beta} \gamma^{\beta} \right>\right>_{0, s} \\
 && - \sum_{\alpha, \beta} (b_{\alpha}+b_{\nu} -2 b_{1} -1) 
        \left<\left< \gamma_{\mu} E \gamma^{\alpha} \right>\right>_{0, s}
        \left<\left< \gamma_{\alpha} \gamma_{\beta}
                \gamma^{\beta} \gamma_{\nu} \right>\right>_{0, s} \\
 && - \sum_{\alpha, \beta} (b_{\alpha}+b_{\mu} - 2 b_{1} -1) 
        \left<\left< \gamma_{\nu} E \gamma^{\alpha} \right>\right>_{0, s}
        \left<\left< \gamma_{\alpha} \gamma_{\beta}
                \gamma^{\beta} \gamma_{\mu} \right>\right>_{0, s} \\
 && + \sum_{\alpha, \beta} 2(2b_{\beta}+b_{\mu}+b_{\nu} - 2 b_{1} -2) 
        \left<\left< E \gamma^{\beta} \gamma^{\alpha} \right>\right>_{0, s}
        \left<\left< \gamma_{\alpha} \gamma_{\beta}
                \gamma_{\mu} \gamma_{\nu} \right>\right>_{0, s} \\
 && + \sum_{\alpha, \beta} 2(b_{\alpha} - b_{1})
                        (b_{\mu} + b_{\nu} - b_{\alpha} - b_{1}) 
        \left<\left< \gamma_{\mu} \gamma_{\nu} \gamma^{\alpha} 
                                \right>\right>_{0, s}
        \left<\left< \gamma_{\alpha} \gamma_{\beta}
                \gamma^{\beta} \right>\right>_{0, s} \\
 && - \sum_{\alpha, \beta} 
        4(b_{\mu}-b_{\alpha}-b_{\beta} - b_{1} + 1)
        (b_{\nu} + b_{\alpha} + b_{\beta} - b_{1} -1) \\
 && \hspace{40pt}
        \left<\left< \gamma_{\mu}  \gamma^{\alpha} \gamma^{\beta}
                        \right>\right>_{0, s}
        \left<\left< \gamma_{\alpha} \gamma_{\beta}
                 \gamma_{\nu} \right>\right>_{0, s}. 
\end{eqnarray*}
Applying Lemma~\ref{lem:derWDVV} (i)  with $u = \gamma_{\nu}$, 
$v = E$, $w_{1} = \gamma_{\mu}$, 
$w_{2} = \gamma_{\beta}$ and  
$w_{3} = \gamma^{\beta}$ to the term
\[     \left<\left< E \gamma^{\beta} \gamma^{\alpha} \right>\right>_{0, s}
        \left<\left< \gamma_{\alpha} \gamma_{\beta}
                \gamma_{\mu} \gamma_{\nu} \right>\right>_{0, s}, \]
and using the symmetry of this expression with respect to $\gamma_{\mu}$
and $\gamma_{\nu}$, then using
Lemma~\ref{lem:sEulerCorr} to 4-point functions which
involve only one $E$ and simplifying, we obtain the desired formula.
$\Box$

We can simplify the formula in Lemma~\ref{lem:G0munu} by the following
simple observation: 
\begin{lem} \label{lem:lowerupper}
For any vector fields $v_{1}, \ldots, v_{k}$ on
the small phase space,
\[ \sum_{\alpha} b_{\alpha} 
	\left< \left< \gamma_{\alpha} \gamma^{\alpha} v_{1} \cdots v_{k}
	\right>\right>_{g, s}
= \frac{1}{2} \sum_{\alpha}  
	\left< \left< \gamma_{\alpha} \gamma^{\alpha} v_{1} \cdots v_{k}
	\right>\right>_{g, s}. \]
\end{lem}
{\bf Proof}:
Since for any $\alpha$ and $\beta$,
$b_{\alpha} \eta^{\alpha \beta} \neq 0$ implies $b_{\alpha} = 1-b_{\beta}$,
we have
\begin{eqnarray*}
 \sum_{\alpha} b_{\alpha} 
	\left< \left< \gamma_{\alpha} \gamma^{\alpha} v_{1} \cdots v_{k}
	\right>\right>_{g, s} 
& = & \sum_{\alpha, \beta} b_{\alpha} \eta^{\alpha \beta} 
	\left< \left< \gamma_{\alpha} \gamma_{\beta} v_{1} \cdots v_{k}
	\right>\right>_{g, s} \\
&=& \sum_{\alpha, \beta} (1-b_{\beta}) \eta^{\alpha \beta} 
	\left< \left< \gamma_{\alpha} \gamma_{\beta} v_{1} \cdots v_{k}
	\right>\right>_{g, s} \\
&=& \sum_{\beta} (1-b_{\beta})  
	\left< \left< \gamma^{\beta} \gamma_{\beta} v_{1} \cdots v_{k}
	\right>\right>_{g, s}.
\end{eqnarray*}
The lemma follows.
$\Box$

Since
\[ G_{0}(E^{m}, E^{k}, E, E) = \sum_{\mu, \nu}
	x_{m}^{\mu} x_{k}^{\nu} G_{0}(\gamma_{\mu}, \gamma_{\nu}, E, E), \]
where $x_{m}^{\mu}$ is defined by (\ref{eqn:EularCoef}),
an immediate consequence of Lemma~\ref{lem:G0munu} and 
Lemma~\ref{lem:lowerupper} 
is the following
\begin{lem} \label{lem:G0mk}
\begin{eqnarray*}
G_{0}(E^{m}, E^{k}, E, E)
 &=& \sum_{\beta} \left<\left< E E E^{m+k} 
                \gamma_{\beta} \gamma^{\beta} \right>\right>_{0, s} \\
 && + \sum_{\alpha, \beta, \mu} (b_{\mu}-b_{\alpha}) 
		x_{m}^{\mu}
        \left<\left< \gamma_{\mu} E \gamma^{\alpha} \right>\right>_{0, s}
        \left<\left< \gamma_{\alpha} \gamma_{\beta}
                \gamma^{\beta} E^{k} \right>\right>_{0, s} \\
 && + \sum_{\alpha, \beta, \nu} (b_{\nu}-b_{\alpha}) 
		x_{k}^{\nu}
        \left<\left< \gamma_{\nu} E \gamma^{\alpha} \right>\right>_{0, s}
        \left<\left< \gamma_{\alpha} \gamma_{\beta}
                \gamma^{\beta} E^{m} \right>\right>_{0, s} \\
 && + \sum_{\alpha, \beta} (b_{\alpha} - b_{1})
                        (2-2b_{\alpha}+2b_{1}) 
        x_{m+k}^{\alpha}
        \left<\left< \gamma_{\alpha} \gamma_{\beta}
                \gamma^{\beta} \right>\right>_{0, s} \\
 && + \sum_{\beta} \left(2b_{1} - 6 + 12 b_{\beta}^{2}\right)
	\left<\left< E^{m+k} \gamma_{\beta}
                \gamma^{\beta} \right>\right>_{0, s} \\
 && + \sum_{\alpha, \beta} 12b_{\alpha}b_{\beta}
	 \left<\left< E^{m} \gamma^{\alpha}
                \gamma^{\beta} \right>\right>_{0, s} 
	\left<\left< \gamma_{\alpha}
                \gamma_{\beta} E^{k} \right>\right>_{0, s} \\
 && + \sum_{\beta, \mu} (b_{\mu}^{2} - b_{\mu}) x_{m}^{\mu}
	\left<\left< \gamma_{\mu} E^{k} \left( \gamma_{\beta}
               \bullet \gamma^{\beta} \right)\right>\right>_{0, s} \\
 && + \sum_{\beta, \nu} (b_{\nu}^{2} - b_{\nu}) x_{k}^{\nu}
	\left<\left< \gamma_{\nu} E^{m} \left( \gamma_{\beta}
               \bullet \gamma^{\beta} \right)\right>\right>_{0, s} \\
 && - \sum_{\beta, \mu, \nu} 2b_{\mu}b_{\nu} 
	x_{m}^{\mu} x_{k}^{\nu}
        \left<\left< \gamma_{\mu} \gamma_{\nu}
	\left( \gamma_{\beta} \bullet \gamma^{\beta} \right)
                        \right>\right>_{0, s}.
\end{eqnarray*}
\end{lem}

To simplify this formula, we need to compute 
$\sum_{\beta}\left<\left< \gamma_{\alpha} \gamma_{\beta}
                \gamma^{\beta} E^{k} \right>\right>_{0, s}$.
First, we have
\begin{lem} \label{lem:4ptPower}
For any vector field $v$ on the small phase space, let $v^{k}$ be
the $k$-th quantum power of $v$. Then for any $\alpha$, $\beta$,
and $\mu$,
\begin{eqnarray*}
\left<\left< v^{k}\gamma_{\alpha} \gamma_{\beta}
                \gamma_{\mu}  \right>\right>_{0, s}
& = & - \sum_{i=1}^{k-1}  
	\left<\left< v^{k-i} 
		\left( \gamma_{\alpha} \bullet \gamma_{\beta}\bullet 
			v^{i-1} \right) v 
                \gamma_{\mu}  \right>\right>_{0, s} \\
&& + \sum_{i=1}^{k}  
	\left<\left< \left( v^{k-i} \bullet \gamma_{\alpha} \right) 
		\left(\gamma_{\beta}\bullet v^{i-1} \right) v 
                \gamma_{\mu}  \right>\right>_{0, s}.
\end{eqnarray*}
\end{lem}
{\bf Proof}: Since 
\[ \left<\left< v^{k}\gamma_{\alpha} \gamma_{\beta}
                \gamma_{\mu}  \right>\right>_{0, s}
 = \sum_{\sigma} \left<\left< v^{k-1} v \gamma_{\sigma} 
			\right>\right>_{0, s}
	\left<\left< \gamma^{\sigma} \gamma_{\alpha} \gamma_{\beta}
                \gamma_{\mu}  \right>\right>_{0, s},
\]
using Lemma~\ref{lem:derWDVV} (i) to exchange positions of
$v$ and $\gamma_{\alpha}$,
we obtain
\begin{eqnarray*}
\left<\left< v^{k}\gamma_{\alpha} \gamma_{\beta}
                \gamma_{\mu}  \right>\right>_{0, s}
& = & - \left<\left< v^{k-1} 
		\left( \gamma_{\alpha} \bullet \gamma_{\beta} \right) v 
                \gamma_{\mu}  \right>\right>_{0, s} \\
&& + \left<\left< \left( v^{k-1} \bullet \gamma_{\alpha} \right) 
		\gamma_{\beta}  v 
                \gamma_{\mu}  \right>\right>_{0, s} \\
&& + \left<\left< v^{k-1} \gamma_{\alpha} 
		\left(\gamma_{\beta} \bullet v \right)
                \gamma_{\mu}  \right>\right>_{0, s}.
\end{eqnarray*}
The lemma follows by 
repeatedly applying this formula to the last term to decrease the
power of the first $v$ and increase the power of the second $v$.
$\Box$

In the special case when $v = E$, Lemma~\ref{lem:4ptPower} implies
\begin{lem} \label{lem:4ptEk}
For any $\mu$ and $k \geq 1$,
\begin{eqnarray*}
\sum_{\beta} \left<\left< E^{k} \gamma_{\beta} \gamma^{\beta} 
	\gamma_{\mu} \right>\right>_{0, s}
&=& \sum_{\beta} (b_{\mu} - b_{1} - 1 + k) 
	\left<\left<  E^{k-1} \gamma_{\mu} \left(\gamma_{\beta}
		\bullet \gamma^{\beta} \right)\right>\right>_{0, s} \\
&& - \sum_{i=1}^{k-1} \sum_{\sigma, \beta} b_{\sigma} x_{k-i}^{\sigma}
	\left<\left<  E^{i-1} 
		\left( \gamma_{\sigma} \bullet \gamma_{\mu} \right)
	 	\left(\gamma_{\beta}
		\bullet \gamma^{\beta} \right)\right>\right>_{0, s} \\
&& - \sum_{i=1}^{k-1} \sum_{\sigma, \beta} b_{\sigma} 
	\left<\left< E^{k-i} \gamma_{\mu} \gamma_{\sigma}
		  \right>\right>_{0, s}
	\left<\left< \gamma^{\sigma} E^{i-1} 
	 	\left(\gamma_{\beta}
		\bullet \gamma^{\beta} \right)\right>\right>_{0, s}.
\end{eqnarray*}
\end{lem}
{\bf Proof}: We first Apply Lemma~\ref{lem:4ptPower} to 
$\left<\left< E^{k} \gamma_{\beta} \gamma^{\beta} 
	\gamma_{\mu} \right>\right>_{0, s}$, then use
Lemma~\ref{lem:sEulerCorr} to remove $E$ from 4-point functions
in the expressions
\[ \left<\left< E^{k-i} 
                \left( \gamma_{\beta} \bullet \gamma^{\beta}\bullet 
                        E^{i-1} \right) E
                \gamma_{\mu}  \right>\right>_{0, s}
 = \sum_{\sigma, \rho} x_{k-i}^{\sigma} 
	\left<\left< E^{i-1}
                \left( \gamma_{\beta} \bullet \gamma^{\beta}
                         \right) 
                \gamma^{\rho}  \right>\right>_{0, s}
	\left<\left< \gamma_{\sigma} \gamma_{\rho} E
                \gamma_{\mu}  \right>\right>_{0, s} \]
and 
\[ \left<\left< \left( E^{k-i} \bullet \gamma^{\beta} \right)
                \left( \gamma_{\beta} \bullet 
                        E^{i-1} \right) E
                \gamma_{\mu}  \right>\right>_{0, s}
 = \sum_{\sigma, \rho} \left<\left< E^{k-i} \gamma^{\beta}
		\gamma^{\sigma} \right>\right>_{0, s}
	\left<\left< E^{i-1}
                \gamma_{\beta} \gamma^{\rho}  \right>\right>_{0, s}
	\left<\left< \gamma_{\sigma} \gamma_{\rho} E
                \gamma_{\mu}  \right>\right>_{0, s}. \]
The lemma is then obtained by using the fact
\[ \left<\left< \left( v_{1} \bullet v_{2} \right) 
		v_{3} v_{4} \right>\right>_{0, s}
	= \left<\left< v_{1} \left( v_{2} \bullet v_{3} \right) 
		 v_{4} \right>\right>_{0, s} \]
for any vector fields $v_{1}, \ldots, v_{4}$ on the small phase space,
and the fact
\[\sum_{\beta} b_{\beta} \left( \gamma_{\beta} \bullet \gamma^{\beta} \right)
= \frac{1}{2} \sum_{\beta}
	 \left( \gamma_{\beta} \bullet \gamma^{\beta} \right), \]
which follows from Lemma~\ref{lem:lowerupper}.
$\Box$

We also need the following
\begin{lem} \label{lem:5ptEulerPower}
\begin{eqnarray*}
\sum_{\beta} \left<\left< E^{k} E E \gamma_{\beta} \gamma^{\beta} 
	 \right>\right>_{0, s}
&=& \sum_{\beta, \mu} (b_{\mu} - b_{1})(b_{\mu} - 2b_{1} -1)
	x_{k}^{\mu}  \left<\left< \gamma_{\mu} \gamma_{\beta}
		\gamma^{\beta} \right>\right>_{0, s} \\
&& + \sum_{\beta, \mu} b_{\mu}(b_{\mu} - b_{1} +k-1) x_{1}^{\mu}
	\left<\left<  E^{k-1} \gamma_{\mu} \left(\gamma_{\beta}
		\bullet \gamma^{\beta} \right)\right>\right>_{0, s} \\
&& - \sum_{i=1}^{k-1} \sum_{\mu, \nu} b_{\mu} b_{\nu} 
		x_{k-i}^{\mu} x_{1}^{\nu}
	\left<\left<  E^{i-1} 
		\left( \gamma_{\mu} \bullet \gamma_{\nu} \right)
	 	\left(\gamma_{\beta}
		\bullet \gamma^{\beta} \right)\right>\right>_{0, s} \\
&& - \sum_{i=1}^{k-1} \sum_{\mu, \nu} b_{\mu} b_{\nu} x_{1}^{\nu} 
	\left<\left< E^{k-i} \gamma_{\mu} \gamma_{\nu}
		  \right>\right>_{0, s}
	\left<\left< \gamma^{\mu} E^{i-1} 
	 	\left(\gamma_{\beta}
		\bullet \gamma^{\beta} \right)\right>\right>_{0, s}.
\end{eqnarray*}
\end{lem}
{\bf Proof}: Since
\[ \sum_{\beta} \left<\left< E^{k} E E \gamma_{\beta} \gamma^{\beta} 
	 \right>\right>_{0, s}
	= \sum_{\beta, \mu, \nu}
	x_{k}^{\mu} x_{1}^{\nu} 
	 \left<\left< \gamma_{\mu} \gamma_{\nu}
		 E \gamma_{\beta} \gamma^{\beta} 
	 \right>\right>_{0, s}, \]
using Lemma~\ref{lem:sEulerCorr} to remove $E$ in the 5-point
function, we obtain
\[ \sum_{\beta} \left<\left< E^{k} E E \gamma_{\beta} \gamma^{\beta} 
	 \right>\right>_{0, s} = 
	\sum_{\beta, \mu} x_{k}^{\mu} (b_{\mu} - 2 b_{1} -1 ) 
	\left<\left< \gamma_{\mu}
		 E \gamma_{\beta} \gamma^{\beta} 
	 \right>\right>_{0, s}
	+ \sum_{\beta, \nu} x_{1}^{\nu} b_{\nu} 
	\left<\left< E^{k} \gamma_{\nu}
		  \gamma_{\beta} \gamma^{\beta} 
	 \right>\right>_{0, s}.
 \]
Using Lemma~\ref{lem:sEulerCorr} again to remove $E$ in the 4-point
function in the first term and applying Lemma~\ref{lem:4ptEk}
to the second term, we obtain the desired formula.
$\Box$

Now we are ready to prove the following
\begin{thm} \label{thm:phik}
 For any manifold $V$,
\begin{eqnarray*}
\phi_{m} & = & - \frac{1}{24}
		\sum_{k=0}^{m-1} \sum_{\alpha, \beta, \sigma} b_{\alpha}
		\left<\left< \gamma_{1} E^{k} 
			\gamma^{\alpha} \right>\right>_{0, s}
		\left<\left< \gamma_{\alpha} E^{m-1-k} 
			\gamma^{\beta} \right>\right>_{0, s} 
		\left<\left< \gamma_{\beta} \gamma_{\sigma}
			\gamma^{\sigma} \right>\right>_{0, s}  \\
	&& - \frac{1}{4}
		\sum_{k=0}^{m-1} \sum_{\alpha, \beta} b_{\alpha} b_{\beta}
		\left<\left< \gamma_{\alpha} E^{k} 
			\gamma^{\beta} \right>\right>_{0, s}
		\left<\left< \gamma_{\beta} E^{m-1-k} 
			\gamma^{\alpha} \right>\right>_{0, s} \\
	&& + \frac{m}{12} \sum_{\sigma}
		\left<\left< \gamma_{\sigma} E^{m-1}
			\gamma^{\sigma} \right>\right>_{0, s}.
\end{eqnarray*}
\end{thm}
{\bf Proof}:
Using formula (\ref{eqn:derCorr}), we obtain
\begin{eqnarray*} E^{m-1} \phi_{2}
&  =  & - \frac{1}{24}
	\sum_{\beta} \left<\left< E^{m-1} E E \gamma_{\beta} \gamma^{\beta} 
	 \right>\right>_{0, s} 
	- \frac{1}{12} \sum_{\beta, \mu} 
	(b_{1} + 1 - b_{\mu}) x_{m-1}^{\mu}
	\left<\left< \gamma_{\mu} E \gamma_{\beta} \gamma^{\beta} 
	 \right>\right>_{0, s}  \\
&& 	+ \frac{1}{2} \sum_{\beta} 
	\left\{ b_{\beta} (1-b_{\beta}) - \frac{b_{1}+1}{6} \right\}
	\left<\left< E^{m-1} \gamma_{\beta} \gamma^{\beta} 
	 \right>\right>_{0, s}.
\end{eqnarray*} 
We can use Lemma~\ref{lem:sEulerCorr} to remove $E$ in the second term,
then plugging the result and the formula in Lemma~\ref{lem:G0mk}
into the definition of $\phi_{m}$. To simplify the resulting expression
of $\phi_{m}$, we first use
Lemma~\ref{lem:5ptEulerPower} to compute the 5-point functions. We can then
use Lemma~\ref{lem:4ptEk} to compute 4-point functions and obtain
the following formula
\begin{eqnarray*}
&& \sum_{k=1}^{m-2} \sum_{\alpha, \beta, \mu} (b_{\mu} - b_{\alpha})
	x_{m-1-k}^{\mu}
	\left<\left< \gamma_{\mu} E \gamma^{\alpha} 
	 \right>\right>_{0, s}
	\left<\left< \gamma_{\alpha} E^{k} \gamma_{\beta} \gamma^{\beta} 
	 \right>\right>_{0, s} \\
& = & \sum_{k=1}^{m-2} \sum_{\alpha, \beta, \mu} 
	b_{\mu}(b_{\alpha}-1)
	x_{m-1-k}^{\mu}
	\left<\left< \gamma_{\mu} E \gamma^{\alpha} 
	 \right>\right>_{0, s}
	\left<\left< \gamma_{\alpha} E^{k-1} 
	\left(\gamma_{\beta} \bullet \gamma^{\beta} \right)
	 \right>\right>_{0, s} \\
&& + \sum_{k=1}^{m-3} \sum_{\alpha, \beta, \mu}
	b_{\alpha}b_{\mu} \left\{
	x_{m-1-k}^{\alpha}x_{k}^{\mu}
	\left<\left< \gamma_{\alpha} \gamma_{\mu} 
	\left(\gamma_{\beta} \bullet \gamma^{\beta} \right)
	 \right>\right>_{0, s}  \right. \\
&& \hspace{100pt} \left.
	- x_{1}^{\alpha}x_{k}^{\mu}
	\left<\left< \left(\gamma_{\alpha} \bullet \gamma_{\mu} \right)
	E^{m-2-k} 
	\left(\gamma_{\beta} \bullet \gamma^{\beta} \right)
	 \right>\right>_{0, s} \right\} \\  
&& + \sum_{k=1}^{m-3} \sum_{\alpha, \beta, \mu}
	b_{\alpha}b_{\mu} \left\{
	x_{m-1-k}^{\alpha}
	\left<\left< \gamma_{\alpha} E \gamma_{\mu} 
	 \right>\right>_{0, s} 
	\left<\left< \gamma^{\mu} E^{k-1} 
	\left(\gamma_{\beta} \bullet \gamma^{\beta} \right)
	 \right>\right>_{0, s}  \right. \\
&& \hspace{100pt} \left.
	- x_{1}^{\alpha}
	\left<\left< \gamma_{\alpha} E^{m-1-k} \gamma_{\mu} 
	 \right>\right>_{0, s} 
	\left<\left< \gamma^{\mu} E^{k-1} 
	\left(\gamma_{\beta} \bullet \gamma^{\beta} \right)
	 \right>\right>_{0, s} \right\} \\
&& - \sum_{k=1}^{m-2}\sum_{\alpha, \beta} b_{\alpha}^{2} x_{m-k}^{\alpha} 
        \left<\left< \gamma_{\alpha} E^{k-1} \left(\gamma_{\beta}
                \bullet \gamma^{\beta} \right)\right>\right>_{0, s} \\
&& + \sum_{\alpha, \beta} b_{1}b_{\alpha} x_{m-1}^{\alpha}
	\left<\left< \gamma_{\alpha} 
	\gamma_{\beta}  \gamma^{\beta}
	 \right>\right>_{0, s} 
 + \sum_{\alpha, \beta} (m-2-b_{1})b_{\alpha} x_{1}^{\alpha}
	\left<\left< \gamma_{\alpha}E^{m-2} 
	\left(\gamma_{\beta} \bullet \gamma^{\beta} \right)
	 \right>\right>_{0, s}.
\end{eqnarray*}
In this way, we obtain
an expression for $\phi_{m}$ which contains only 3-point functions. 
There are many cancellations
among different terms in this expression. After simplifying it, we obtain
the desired formula.
$\Box$

{\bf Remark}: The expression for $\phi_{m}$ in Theorem~\ref{thm:phik}
is the same as that for $\left<\left< E^{k} \right>\right>_{1, s}$
obtained in \cite{DZ2} (4.42) for the case where the quantum cohomology
of $V$ is semisimple.

\section{A necessary and sufficient condition for the genus-1
	Virasoro conjecture}
\label{sec:g1VirIff}

The main purpose of this section is to prove Theorem~\ref{thm:phi2g1Vir}.
In this section, we will use $\{u^{1}, \ldots, u^{N}\}$ to denote the
coordinate on the small phase space in order to distinguish
the one on the big phase space. In this coordinate,
the vector field $\frac{\partial}{\partial u^{\alpha}}$ is identified
with $\gamma_{\alpha}$. Let 
$u_{\alpha} = \sum_{\beta} \eta_{\alpha \beta} u^{\beta}$. Then
$\frac{\partial}{\partial u_{\alpha}}$ is identified with $\gamma^{\alpha}$.
Let $M$ be an $N \times N$ matrix whose entries are $u_{\alpha \beta}$.
Temporarily, we think of each $u_{\alpha \beta}$ as an independent
variable.
Define 
\[ {\cal F}_{1}(u_{1}, \ldots, u_{N}; M)
	:= \left< e^{\sum_{\alpha} 
		u_{\alpha} \gamma^{\alpha}} \right>_{1}
        + \frac{1}{24} \log \det \left( \eta^{-1}M \right). \]
Then
\begin{equation} \label{eqn:transDer}
\frac{\partial {\cal F}_{1}}{\partial u_{\alpha}}
= \left<\left< \gamma^{\alpha}  \right>\right>_{0, s}
\, \, \, \, \, {\rm and} \, \, \, \, \,
\frac{\partial {\cal F}_{1}}{\partial u_{\alpha \beta}}
= \frac{1}{24}\left(M^{-1}\right)_{\alpha \beta}.
\end{equation}
The genus-1 constitutive relation says that 
${\cal F}_{1}$ is equal to $F_{1}$ after the transformation
\begin{equation} \label{eqn:trans}
	u_{\alpha} =  \left<\left< \gamma_{1} \gamma_{\alpha} 
                        \right>\right>_{0}
	\, \, \, \, \, {\rm and} \, \, \, \, \,
	u_{\alpha \beta} =  \left<\left< \gamma_{1} \gamma_{\alpha} 
                       \gamma_{\beta} \right>\right>_{0}.
\end{equation}
Taking derivative of the genus-1 constitutive relation once, we obtain
\begin{equation} \label{eqn:derCons}
\left<\left< \tau_{m, \alpha} \right>\right>_{1}
= \sum_{\sigma} \left<\left< \gamma_{1}\tau_{m, \alpha} 
                       \gamma_{\sigma} \right>\right>_{0}
	\frac{\partial {\cal F}_{1}}{\partial u_{\sigma}}
	+ \sum_{\sigma, \rho}\left<\left< \gamma_{1}\tau_{m, \alpha} 
       \gamma_{\sigma} \gamma_{\rho} \right>\right>_{0} 
              \frac{\partial {\cal F}_{1}}{\partial u_{\sigma \rho}}
\end{equation}
for any $m$ and $\alpha$.
On the other hand, the genus-0 constitutive relation says, in particular,
that
\begin{equation} \label{eqn:trans2pt}
\left. \left<\left< \gamma_{\alpha}
                \gamma_{\beta}  \right>\right>_{0, s} 
		\right|_{u_{\sigma} =  \left<\left< \gamma_{1} \gamma_{\sigma} 
                        \right>\right>_{0}}
= \left<\left< \gamma_{\alpha}
                \gamma_{\beta}  \right>\right>_{0}. 
\end{equation}
Taking derivative of this relation once, we get
\begin{equation} \label{eqn:trans3pt}
\left. \left<\left< \gamma_{\alpha}
                \gamma_{\beta} \gamma_{\mu} \right>\right>_{0, s} 
		\right|_{u_{\sigma} =  \left<\left< \gamma_{1} \gamma_{\sigma} 
                        \right>\right>_{0}}
= \sum_{\nu} \left. \left(M^{-1} \eta \right)_{\nu \mu} 
		\right|_{u_{\sigma \rho} = 
		\left<\left< \gamma_{1} \gamma_{\sigma} 
                       \gamma_{\rho} \right>\right>_{0}}
	\left<\left< \gamma_{\alpha}
                \gamma_{\beta} \gamma_{\nu} \right>\right>_{0}. 
\end{equation}
Moreover combining equation (\ref{eqn:trans2pt}) with 
Lemma~\ref{lem:EulerCorr} (iii) and Lemma~\ref{lem:sEulerCorr} (iii), we obtain
\begin{equation} \label{eqn:transEuler}
\left. \left<\left< \gamma_{\alpha} E
                \gamma_{\beta}  \right>\right>_{0, s} 
		\right|_{u_{\sigma} =  \left<\left< \gamma_{1} \gamma_{\sigma} 
                        \right>\right>_{0}}
= \left<\left< \gamma_{\alpha} {\cal X}
                \gamma_{\beta}  \right>\right>_{0}. 
\end{equation}
The following lemma will be useful in the proof of Theorem~\ref{thm:phi2g1Vir}.
\begin{lem} \label{lem:gamma1Euler}
\begin{eqnarray*}
&& \sum_{\alpha, \beta, \mu_{1}, \ldots, \mu_{k-1}}
	\left( M^{-1} \right)_{\alpha \beta} \gamma_{1} \left\{ 
	\left<\left<  \gamma_{\alpha} {\cal X} \gamma^{\mu_{1}} 
		\right>\right>_{0} 
	\left<\left<  \gamma_{\mu_{1}} {\cal X} \gamma^{\mu_{2}} 
		\right>\right>_{0} 
	\cdots \left<\left<  \gamma_{\mu_{k-1}} {\cal X} \gamma_{\beta} 
		\right>\right>_{0} 
	\right\}  \\
& = & k \sum_{\mu_{1}, \ldots, \mu_{k-1}}  
	\left<\left<  \gamma_{\mu_{k-1}} {\cal X} \gamma^{\mu_{1}} 
		\right>\right>_{0} 
	\left<\left<  \gamma_{\mu_{1}} {\cal X} \gamma^{\mu_{2}} 
		\right>\right>_{0} 
	\cdots \left<\left<  \gamma_{\mu_{k-2}} {\cal X} \gamma^{\mu_{k-1}} 
		\right>\right>_{0}.
\end{eqnarray*}
\end{lem}
{\bf Proof}: By Lemma~\ref{lem:EulerCorr} (iii),
\[ \gamma_{1} \left<\left<  \gamma_{\alpha} {\cal X} \gamma_{\beta} 
		\right>\right>_{0} 
	= (b_{\alpha} + b_{\beta})
	\left<\left< \gamma_{1} \gamma_{\alpha} \gamma_{\beta} 
		\right>\right>_{0} 
	= (b_{\alpha} + b_{\beta}) u_{\alpha \beta}. \]
Therefore
\begin{eqnarray*}
&& \sum_{\alpha, \beta, \mu_{1}, \ldots, \mu_{k-1}}
	\left( M^{-1} \right)_{\alpha \beta} \gamma_{1} \left\{ 
	\left<\left<  \gamma_{\alpha} {\cal X} \gamma^{\mu_{1}} 
		\right>\right>_{0} 
	\left<\left<  \gamma_{\mu_{1}} {\cal X} \gamma^{\mu_{2}} 
		\right>\right>_{0} 
	\cdots \left<\left<  \gamma_{\mu_{k-1}} {\cal X} \gamma_{\beta} 
		\right>\right>_{0} 
	\right\}  \\
& = & 2k \sum_{\mu_{1}, \ldots, \mu_{k-1}}  b_{\mu_{1}} 
	\left<\left<  \gamma_{\mu_{k-1}} {\cal X} \gamma^{\mu_{1}} 
		\right>\right>_{0} 
	\left<\left<  \gamma_{\mu_{1}} {\cal X} \gamma^{\mu_{2}} 
		\right>\right>_{0} 
	\cdots \left<\left<  \gamma_{\mu_{k-2}} {\cal X} \gamma^{\mu_{k-1}} 
		\right>\right>_{0}.
\end{eqnarray*}
In this calculation, one needs to switch the position of $\gamma_{1}$ and
that ot ${\cal X}$ by using the generalized WDVV equation 
so that $\gamma_{1}$ can be pushed to the beginning or the end of
the chain of the multiplications of 3-point functions. In this way we can 
always create entries of $M$ which can be used to eliminate entries
of $M^{-1}$. Moreover, by interchanging all upper indices with
the corresponding lower indices, we obtain
\begin{eqnarray*}
&& \sum_{\mu_{1}, \ldots, \mu_{k-1}}  b_{\mu_{1}} 
	\left<\left<  \gamma_{\mu_{k-1}} {\cal X} \gamma^{\mu_{1}} 
		\right>\right>_{0} 
	\left<\left<  \gamma_{\mu_{1}} {\cal X} \gamma^{\mu_{2}} 
		\right>\right>_{0} 
	\cdots \left<\left<  \gamma_{\mu_{k-2}} {\cal X} \gamma^{\mu_{k-1}} 
		\right>\right>_{0} \\
&=& \frac{1}{2}
	\sum_{\mu_{1}, \ldots, \mu_{k-1}} 
	\left<\left<  \gamma_{\mu_{k-1}} {\cal X} \gamma^{\mu_{1}} 
		\right>\right>_{0} 
	\left<\left<  \gamma_{\mu_{1}} {\cal X} \gamma^{\mu_{2}} 
		\right>\right>_{0} 
	\cdots \left<\left<  \gamma_{\mu_{k-2}} {\cal X} \gamma^{\mu_{k-1}} 
		\right>\right>_{0}.
\end{eqnarray*}
The lemma then follows. $\Box$

Recall that ${\cal L}_{1}$ is the vector field on the big phase space
which is defined to be the first derivative part of the $L_{1}$ operator.
The genus-1 $L_{1}$ constraint is $\Psi_{1, 1} = 0$, where
\[ 
 \Psi_{1, 1} = \left<\left< {\cal L}_{1} \right>\right>_{1}
	+ \frac{1}{2}\sum_{\alpha} b_{\alpha} (1-b_{\alpha}) 
	\left\{\left<\left< \gamma_{\alpha} \gamma^{\alpha} \right>\right>_{0}
	+ 2  \left<\left< \gamma^{\alpha} \right>\right>_{0}
		 \left<\left< \gamma_{\alpha} \right>\right>_{1} 
		\right\}.
\]
We have the following
\begin{pro} \label{pro:g1L1}
\begin{eqnarray*}
 \Psi_{1, 1}& = & \left. \left\{ - \left<\left< E^{2} \right>\right>_{1, s}
	\, + \, \phi_{2}
	\right\} 
	\right|_{u_{\sigma} =  \left<\left< \gamma_{1} \gamma_{\sigma} 
                        \right>\right>_{0}}.
\end{eqnarray*}
\end{pro}
{\bf Proof}: Applying equation (\ref{eqn:transDer}) and
(\ref{eqn:derCons}) to each genus-1 1-point
function in $\Psi_{1, 1}$, we obtain
\begin{eqnarray} 
 \Psi_{1, 1}& = &\sum_{\sigma} 
		\left\{ \left<\left< \gamma_{1} {\cal L}_{1} 
                       \gamma_{\sigma} \right>\right>_{0}
		+ \sum_{\alpha} b_{\alpha} (1-b_{\alpha})
		\left<\left< \gamma^{\alpha} \right>\right>_{0}
		\left<\left< \gamma_{1} \gamma_{\alpha} 
                       \gamma_{\sigma} \right>\right>_{0} \right\}
         \left\{\left. \left<\left< \gamma^{\sigma}  \right>\right>_{0, s}
			\right|_{ u_{\beta} =  
		\left<\left< \gamma_{1} \gamma_{\beta} 
                        \right>\right>_{0}} \right\}
		\nonumber \\
	&& 
	+ \frac{1}{24}\sum_{\sigma, \rho} \left\{
	\left<\left< \gamma_{1} {\cal L}_{1}
       \gamma_{\sigma} \gamma_{\rho} \right>\right>_{0} 
	+\sum_{\alpha} b_{\alpha} (1-b_{\alpha})
		\left<\left< \gamma^{\alpha} \right>\right>_{0}
	\left<\left< \gamma_{1} \gamma_{\alpha}
       \gamma_{\sigma} \gamma_{\rho} \right>\right>_{0}  \right\}
		\left(M^{-1}\right)_{\sigma \rho}
		\nonumber \\
	&& + \frac{1}{2}\sum_{\alpha} b_{\alpha} (1-b_{\alpha}) 
	\left<\left< \gamma_{\alpha} \gamma^{\alpha} \right>\right>_{0},
		\label{eqn:g1L1}
\end{eqnarray}
where the entries of $M$ are $u_{\alpha \beta} = 
	\left<\left< \gamma_{1} \gamma_{\alpha} \gamma_{\beta} 
		\right>\right>_{0}$.
By the second equation of (\ref{eqn:VirVect3pt}) and
equation (\ref{eqn:transEuler}), the first line of the right hand side
is equal to $\left.  - \left<\left< E^{2} \right>\right>_{1, s}
		\right|_{u_{\sigma} =  \left<\left< \gamma_{1} \gamma_{\sigma} 
                        \right>\right>_{0}}$.
Now we compute the second line. Since
\[ \left<\left< \gamma_{1} {\cal L}_{1}
       \gamma_{\sigma} \gamma_{\rho} \right>\right>_{0}
	= \gamma_{1} \left<\left<  {\cal L}_{1}
       \gamma_{\sigma} \gamma_{\rho} \right>\right>_{0}
	- b_{1}(b_{1}+1) \left<\left< \tau_{1, 1}
       \gamma_{\sigma} \gamma_{\rho} \right>\right>_{0}
	- (2b_{1}+1) \sum_{\alpha}{\cal C}_{1}^{\alpha} 
	\left<\left< \gamma_{\alpha}
       \gamma_{\sigma} \gamma_{\rho} \right>\right>_{0},
\]
by Lemma~\ref{lem:gamma1Euler} and
the second equation of (\ref{eqn:VirVect3pt}), 
the second line of (\ref{eqn:g1L1}) is equal to
\begin{eqnarray}
&& - \frac{1}{12} \sum_{\alpha} 
	\left<\left< \gamma_{\alpha} \gamma^{\alpha} \right>\right>_{0}
-\frac{1}{24}\sum_{\sigma, \rho, \alpha}  
	\left(M^{-1}\right)_{\sigma \rho} 
	\left\{ 
	\left( b_{\alpha} (1-b_{\alpha}) 
		+ b_{1} (b_{1}+1) \right)
	\left<\left< \gamma_{1} \gamma_{\alpha} \right>\right>_{0}
	\left<\left< \gamma^{\alpha} \gamma_{\sigma} \gamma_{\rho}
		 \right>\right>_{0} \right.
	\nonumber \\
	&& \hspace{200pt} \left.
 	+ (2b_{1} +1) {\cal C}_{1}^{\alpha} 
	\left<\left< \gamma_{\alpha}
       \gamma_{\sigma} \gamma_{\rho} \right>\right>_{0} 
 	\right\}.    \label{eqn:2ndg1L1}
\end{eqnarray}
On the other hand, by
Lemma~\ref{lem:sEulerCorr} (iv) 
\begin{eqnarray*}
\sum_{\alpha} 
    \left<\left< E E \gamma_{\alpha} \gamma^{\alpha} \right>\right>_{0, s}
	& = & \sum_{\alpha, \beta}
      \left<\left< E \gamma_{1} \gamma^{\beta} \right>\right>_{0, s}
	\left<\left< \gamma_{\beta} E \gamma_{\alpha} 
		\gamma^{\alpha} \right>\right>_{0, s} \\
	&=& \sum_{\alpha, \beta} (b_{\beta} - b_{1})
      \left<\left< E \gamma_{1} \gamma^{\beta} \right>\right>_{0, s}
	\left<\left< \gamma_{\beta} \gamma_{\alpha} 
		\gamma^{\alpha} \right>\right>_{0, s}.
\end{eqnarray*}
By  equation (\ref{eqn:trans3pt}) and (\ref{eqn:transEuler}),
\[
\sum_{\alpha} \left.
         \left<\left< E E \gamma_{\alpha} \gamma^{\alpha} \right>\right>_{0, s}
	\right|_{u_{\sigma} =  \left<\left< \gamma_{1} \gamma_{\sigma} 
                        \right>\right>_{0}}
= \sum_{\beta, \sigma, \rho} (b_{\beta} - b_{1})
      \left<\left< X \gamma_{1} \gamma^{\beta} \right>\right>_{0}
	\left<\left< \gamma_{\beta} \gamma_{\sigma} 
		\gamma_{\rho} \right>\right>_{0} 
	\left(M^{-1}\right)_{\sigma \rho}.
\]
By Lemma~\ref{lem:EulerCorr} (iii), 
\begin{eqnarray*} &&
 \sum_{\alpha} b_{\alpha} \left<\left< \gamma_{1} {\cal X} \gamma^{\alpha}
        \right>\right>_{0}
	 \\
&=& \sum_{\alpha} \left\{  - b_{1} 
	\left<\left< \gamma_{1} {\cal X} \gamma^{\alpha}\right>\right>_{0}
	+(2b_{1} +1) {\cal C}_{1}^{\alpha} 
	+ \left( b_{\alpha} (1-b_{\alpha}) 
		+ b_{1} (b_{1}+1) \right)
	\left<\left< \gamma_{1} \gamma^{\alpha} \right>\right>_{0} \right\}.
\end{eqnarray*}
Moreover
\begin{eqnarray*}
\sum_{\beta, \sigma, \rho}
      \left<\left< X \gamma_{1} \gamma^{\beta} \right>\right>_{0}
	\left<\left< \gamma_{\beta} \gamma_{\sigma} 
		\gamma_{\rho} \right>\right>_{0} 
	\left(M^{-1}\right)_{\sigma \rho}
&=& \sum_{\beta, \sigma, \rho} 
      \left<\left< X \gamma_{\sigma} \gamma^{\beta} \right>\right>_{0}
	\left<\left< \gamma_{\beta} \gamma_{1} 
		\gamma_{\rho} \right>\right>_{0} 
	\left(M^{-1}\right)_{\sigma \rho} \\
&=& \sum_{\sigma} 
      \left<\left< X \gamma_{\sigma} \gamma^{\sigma} \right>\right>_{0}
	\\
&=& \sum_{\sigma} 
      \left<\left< \gamma_{\sigma} \gamma^{\sigma} \right>\right>_{0}.
\end{eqnarray*}
Therefore we have 
\begin{eqnarray*}
&& \sum_{\alpha} \left.
         \left<\left< E E \gamma_{\alpha} \gamma^{\alpha} \right>\right>_{0, s}
	\right|_{u_{\sigma} =  \left<\left< \gamma_{1} \gamma_{\sigma} 
                        \right>\right>_{0}} \\
&=& \sum_{\alpha, \sigma, \rho} \left\{ 
	(2b_{1} +1) {\cal C}_{1}^{\alpha} 
	+ \left( b_{\alpha} (1-b_{\alpha}) 
		+ b_{1} (b_{1}+1) \right)
	\left<\left< \gamma_{1} \gamma^{\alpha} \right>\right>_{0} \right\}
	\left<\left< \gamma_{\alpha} \gamma_{\sigma} 
		\gamma_{\rho} \right>\right>_{0} 
	\left(M^{-1}\right)_{\sigma \rho} \\
&&	- 2 b_{1} \sum_{\sigma} 
	\left<\left< \gamma_{\sigma} \gamma^{\sigma} \right>\right>_{0}.
\end{eqnarray*}
Comparing this equation with (\ref{eqn:2ndg1L1}) and using
(\ref{eqn:trans2pt}), we obtain the desired formula.
$\Box$

We next prove the analogue of this proposition for the genus-1 $L_{2}$
constraint. We need the following
\begin{lem} \label{lem:SomeIdentity}
\begin{eqnarray*}
&{\rm (i)} & \sum_{\alpha, \beta} 
	b_{\alpha}\left<\left< \gamma_{\alpha} \gamma^{\beta} 
		\right>\right>_{0} 
	\left<\left< \gamma_{\beta} \gamma^{\alpha} 
		\right>\right>_{0} 
	= \frac{1}{2}
		\sum_{\alpha, \beta} 
	\left<\left< \gamma_{\alpha} \gamma^{\beta} 
		\right>\right>_{0} 
	\left<\left< \gamma_{\beta} \gamma^{\alpha} 
		\right>\right>_{0},   \\
&{\rm (ii)} & \sum_{\alpha, \beta} 
	b_{\alpha}^{3}\left<\left< \gamma_{\alpha} \gamma^{\beta} 
		\right>\right>_{0} 
	\left<\left< \gamma_{\beta} \gamma^{\alpha} 
		\right>\right>_{0} 
	= \sum_{\alpha, \beta} 
	\left(- \frac{1}{4} + \frac{3}{2} b_{\alpha}^{2}\right)
	\left<\left< \gamma_{\alpha} \gamma^{\beta} 
		\right>\right>_{0} 
	\left<\left< \gamma_{\beta} \gamma^{\alpha} 
		\right>\right>_{0},   \\
&{\rm (iii)}& \sum_{\alpha, \beta} (b_{\alpha})^{k}
	{\cal C}_{\alpha}^{\beta} 
	\left<\left< \gamma_{\beta} \gamma^{\alpha} 
		\right>\right>_{0}  = 0
	\, \, \, \, {\rm if} \, \, \, \, k \, \, \, \,
	{\rm is \, \, \, \, odd}.
\end{eqnarray*}
\end{lem}
{\bf Proof}: Interchanging the upper indices and lower indices
in the expression  \\
$\sum_{\alpha, \beta} 
	b_{\alpha}\left<\left< \gamma_{\alpha} \gamma^{\beta} 
		\right>\right>_{0} 
	\left<\left< \gamma_{\beta} \gamma^{\alpha} 
		\right>\right>_{0}$ and using the fact that
$b_{\alpha} \eta^{\alpha \beta} \neq 0$ implies $b_{\beta} = 1 - b_{\alpha}$, 
we obtain
\[ \sum_{\alpha, \beta} 
	b_{\alpha}\left<\left< \gamma_{\alpha} \gamma^{\beta} 
		\right>\right>_{0} 
	\left<\left< \gamma_{\beta} \gamma^{\alpha} 
		\right>\right>_{0} 
	= \sum_{\alpha, \beta} (1-b_{\alpha})
	\left<\left< \gamma_{\beta} \gamma^{\alpha} 
		\right>\right>_{0} 
	\left<\left< \gamma_{\alpha} \gamma^{\beta} 
		\right>\right>_{0} . \]
This implies (i). Similarly we have
\[ \sum_{\alpha, \beta} 
	b_{\alpha}^{3}\left<\left< \gamma_{\alpha} \gamma^{\beta}
		\right>\right>_{0} 
	\left<\left< \gamma_{\beta} \gamma^{\alpha} 
		\right>\right>_{0} 
	= \sum_{\alpha, \beta} (1-b_{\alpha})^{3}
	\left<\left< \gamma_{\beta} \gamma^{\alpha} 
		\right>\right>_{0} 
	\left<\left< \gamma_{\alpha} \gamma^{\beta} 
		\right>\right>_{0} . \]
Together with (i), this implies (ii).
Using the fact
$b_{\alpha} {\cal C}_{\alpha \beta} \neq 0$ implies 
$b_{\beta} = - b_{\alpha}$, we have
\[\sum_{\alpha, \beta} (b_{\alpha})^{k}
	{\cal C}_{\alpha}^{\beta} 
	\left<\left< \gamma_{\beta} \gamma^{\alpha} 
		\right>\right>_{0}  
 = \sum_{\alpha, \beta} (b_{\alpha})^{k}
	{\cal C}_{\alpha \beta} 
	\left<\left< \gamma^{\beta} \gamma^{\alpha} 
		\right>\right>_{0} 
 = \sum_{\alpha, \beta} (- b_{\beta})^{k}
	{\cal C}_{\alpha \beta} 
	\left<\left< \gamma^{\beta} \gamma^{\alpha} 
		\right>\right>_{0} . \]
Interchanging $\alpha$ with $\beta$, we have
\[\sum_{\alpha, \beta} (b_{\alpha})^{k}
	{\cal C}_{\alpha}^{\beta} 
	\left<\left< \gamma_{\beta} \gamma^{\alpha} 
		\right>\right>_{0}  = 
  (-1)^{k} \sum_{\alpha, \beta} (b_{\alpha})^{k}
	{\cal C}_{\alpha}^{\beta} 
	\left<\left< \gamma_{\beta} \gamma^{\alpha} 
		\right>\right>_{0}. \]
This implies (iii). $\Box$

 The genus-1 $L_{2}$ constraint is the equation 
$\Psi_{1, 2} = 0$ where
\begin{eqnarray*}
\Psi_{1, 2} & = &  \left<\left< {\cal L}_{2} \right>\right>_{1}
 	+ \sum_{\alpha} b_{\alpha} (1-b_{\alpha}^{2}) \left\{ 
        \left<\left< \tau_{1, \alpha} \gamma^{\alpha} \right>\right>_{0}
	+ \left<\left< \tau_{1, \alpha} \right>\right>_{0}
                 \left<\left< \gamma^{\alpha} \right>\right>_{1}
	+ \left<\left< \tau_{1, \alpha} \right>\right>_{1}
                 \left<\left< \gamma^{\alpha} \right>\right>_{0} \right\} \\
	&&
	- \frac{1}{2}\sum_{\alpha, \beta} (3b_{\alpha}^{2} -1)
		{\cal C}_{\alpha}^{\beta} \left\{
		\left<\left< \gamma^{\alpha} \gamma_{\beta} \right>\right>_{0}
		+ 2 \left<\left< \gamma^{\alpha} \right>\right>_{1}
                 \left<\left< \gamma_{\beta} \right>\right>_{0} \right\}.
\end{eqnarray*}
We have the following
\begin{pro} \label{pro:g1L2}
\begin{eqnarray*}
\Psi_{1, 2}& = & \left.\left\{ - \left<\left< E^{3} \right>\right>_{1, s}
		+ \phi_{3}
		\right\}
		\right|_{u_{\sigma} =  \left<\left< \gamma_{1} \gamma_{\sigma} 
                        \right>\right>_{0}}.
\end{eqnarray*}
\end{pro}
{\bf Proof}: 
Applying equation (\ref{eqn:transDer}) and
(\ref{eqn:derCons}) to each genus-1 1-point
function in $\Psi_{1, 2}$, using equation (\ref{eqn:3ptL2}) and the
fact
\begin{eqnarray*}
 \left<\left< \gamma_{1} {\cal L}_{2}
       \gamma_{\sigma} \gamma_{\rho} \right>\right>_{0}
       & = & \gamma_{1} \left<\left<  {\cal L}_{2}
       \gamma_{\sigma} \gamma_{\rho} \right>\right>_{0}
        - b_{1}(b_{1}+1)(b_{1}+2) \left<\left< \tau_{2, 1}
       \gamma_{\sigma} \gamma_{\rho} \right>\right>_{0} \\
	&& - \sum_{\beta} (3b_{1}^{2} + 6 b_{1} +2) {\cal C}_{1}^{\beta}
		\left<\left< \tau_{1, \beta}
       \gamma_{\sigma} \gamma_{\rho} \right>\right>_{0} \\
       && - \sum_{\beta} 3(b_{1}+1) 
	({\cal C}^{2})_{1}^{\beta}  
        \left<\left< \gamma_{\beta}
       \gamma_{\sigma} \gamma_{\rho} \right>\right>_{0},
\end{eqnarray*}
then applying Lemma~\ref{lem:gamma1Euler}, equation (\ref{eqn:transEuler})
and the genus-0 topological recursion relation,
 we obtain
\begin{eqnarray}
\Psi_{1, 2}& = & \left.\left\{ - \left<\left< E^{3} \right>\right>_{1, s}
                - \frac{1}{8} \sum_{\alpha} 
		\left<\left< E^{2} \gamma_{\alpha} \gamma^{\alpha}
			 \right>\right>_{1, s}
                \right\}
                \right|_{u_{\sigma} =  \left<\left< \gamma_{1} \gamma_{\sigma} 
                        \right>\right>_{0}} \nonumber \\
		&& + \frac{1}{24} \sum_{\mu, \nu, \beta}
		\left(M^{-1}\right)_{\mu \nu} 
			\left<\left< \gamma_{\mu} 
                       \gamma_{\nu}\gamma^{\beta}  \right>\right>_{0} 
			\nonumber \\
		&& \hspace{50pt}
			\left\{ 
			 \sum_{\alpha} (3 b_{\alpha}^{2} -1)
			{\cal C}_{\alpha \beta} 
			\left<\left< \gamma_{1} \gamma^{\alpha}
				\right>\right>_{0} \right.
			+ \sum_{\alpha} b_{\alpha} (b_{\alpha}^{2} -1) 
			\left<\left< \gamma_{1} \gamma^{\alpha}
				\right>\right>_{0}
                \left<\left< \gamma_{\alpha} \gamma_{\beta} 
                       \right>\right>_{0}  
			\nonumber \\
		&& \hspace{60pt} 
		+ b_{\beta} (b_{\beta}^{2} -1) 
			\left<\left< \gamma_{1} \tau_{1, \beta}
				\right>\right>_{0}
		 - b_{1}(b_{1}+1)(b_{1}+2) \left<\left< \tau_{1, 1}
       		\gamma_{\beta}  \right>\right>_{0} 
			\nonumber \\
		&& \hspace{60pt} \left.
		 - \sum_{\alpha} (3b_{1}^{2} + 6 b_{1} +2) 
			{\cal C}_{1}^{\alpha}
                	\left<\left< 
       			\gamma_{\alpha} \gamma_{\beta} \right>\right>_{0} 
		 - 3(b_{1}+1) 
        		({\cal C}^{2})_{1 \beta}  
        		\right\}
		\nonumber \\
	&& - \sum_{\beta} b_{\beta} (b_{\beta}^{2} -1) 
			\left<\left< \tau_{1, \beta} \gamma^{\beta} 
                        \right>\right>_{0}  
		- \frac{1}{2}\sum_{\alpha, \beta} (3 b_{\alpha}^{2} -1)
			{\cal C}_{\alpha}^{\beta} 
			\left<\left< \gamma_{\beta} \gamma^{\alpha}
				\right>\right>_{0}.
	\label{eqn:g1L2-1}
\end{eqnarray}
A simple combination of Lemma~\ref{lem:EulerCorr} and the genus-0 topological
recursion relation gives the following (cf. \cite{LT} formula (8) and (9))
\begin{eqnarray*}
 (1 + b_{\alpha} + b_{\beta}) \left<\left< \tau_{1, \alpha} \gamma_{\beta} 
	\right>\right>_{0}  
& = & \sum_{\sigma} \left<\left< \gamma_{\alpha} \gamma^{\sigma} 
	\right>\right>_{0} 
	\left\{ {\cal C}_{\sigma \beta} + (b_{\sigma} + b_{\beta})
	\left<\left< \gamma_{\sigma} \gamma_{\beta} 
                        \right>\right>_{0}  \right\}
	- \sum_{\sigma} {\cal C}_{\alpha}^{\sigma} 
	\left<\left< \gamma_{\sigma} \gamma_{\beta} 
                        \right>\right>_{0}.
\end{eqnarray*} 
This is a special case of the fundamental recursion relation of
\cite{EHX1}. Using this formula, we can express  
2-point correlation functions of type
$\left<\left< \tau_{1, \alpha} \gamma_{\beta} \right>\right>_{0}$ in 
the right hand side of equation (\ref{eqn:g1L2-1}) in terms of correlation
functions only involving $\gamma_{\sigma}$, $\sigma = 1, \ldots, N$.
(In this procedure, first applying Lemma 3.2 in \cite{LT} to shift the level
of descendant in the term
$b_{\beta} (1+ b_{\beta}) 
	\left<\left< \gamma_{1} \tau_{1, \beta} \right>\right>_{0}$ may 
simplify the computation.)
Then a straightforward computation using Lemma~\ref{lem:EulerCorr}
and Lemma~\ref{lem:SomeIdentity} shows that
\begin{eqnarray*}
\Psi_{1, 2}& = & \left.\left\{ - \left<\left< E^{3} \right>\right>_{1, s}
                - \frac{1}{8} \sum_{\alpha} 
                \left<\left< E^{2} \gamma_{\alpha} \gamma^{\alpha}
                         \right>\right>_{1, s}
                \right\}
                \right|_{u_{\sigma} =  \left<\left< \gamma_{1} \gamma_{\sigma} 
                        \right>\right>_{0}} \nonumber \\
                && + \frac{1}{24} \sum_{\mu, \nu, \beta}
                \left(M^{-1}\right)_{\mu \nu} 
                        \left<\left< \gamma_{\mu} 
                       \gamma_{\nu}\gamma^{\beta}  \right>\right>_{0}
			(b_{1} + b_{\alpha} +1 - b_{\beta}) 
			\left<\left< \gamma_{1} {\cal X} 
                       \gamma^{\alpha}  \right>\right>_{0} 
			\left<\left< \gamma_{\alpha} {\cal X} 
                       \gamma_{\beta}  \right>\right>_{0} 
                        \nonumber \\
                && + \sum_{\alpha, \beta}
		\left( \frac{3}{8} - \frac{1}{2} b_{\beta}^{2}
			- \frac{1}{4} b_{\alpha} b_{\beta} \right)
		 \left<\left< \gamma_{\alpha} {\cal X} 
                       \gamma^{\beta}  \right>\right>_{0} 
		\left<\left< \gamma_{\beta} {\cal X} 
                       \gamma^{\alpha}  \right>\right>_{0}. 
\end{eqnarray*}
The proposition then follows from equation (\ref{eqn:trans3pt}),
(\ref{eqn:transEuler}), and Theorem~\ref{thm:phik}. $\Box$

Now we are ready to prove Theorem~\ref{thm:phi2g1Vir}.

{\bf Proof of Theorem~\ref{thm:phi2g1Vir}}:
The string equation implies that the transformation \\
$u^{\alpha} = \left<\left< \gamma_{1} \gamma^{\alpha}  \right>\right>_{0, s}$
is an identity map when the right hand side of this equation is restricted
to the small phase space. Therefore, by Proposition~\ref{pro:g1L1},
the restriction of the genus-1 $L_{1}$ constraint to the small
phase space is equivalent to the condition that
$\left<\left< E^{2}  \right>\right>_{0, s} = \phi_{2}$.
Hence $\left<\left< E^{2}  \right>\right>_{0, s} = \phi_{2}$
 is a necessary condition for the genus-1 Virasoro
conjecture. On the other hand, if 
$\left<\left< E^{2}  \right>\right>_{0, s} = \phi_{2}$, 
Proposition~\ref{pro:g1L1} also implies that the genus-1 $L_{1}$ constraint
is true. Moreover, Theorem~\ref{thm:EulerPhi} and Proposition~\ref{pro:g1L2}
implies that the genus-1 $L_{2}$ constraint is also true. By the virasoro
relation among the $L_{n}$ operators, the genus-1 Virasoro conjecture holds.
$\Box$

\section{Virasoro type relation for $\{ \phi_{k} \}$}
\label{sec:phiVir}

Because of Theorem~\ref{thm:phi2g1Vir}, we are interested in
when the equality $\left<\left<  E^{2} \right>\right>_{1, s} = \phi_{2}$
holds. The Virasoro relation (\ref{eqn:VirEuler}) and 
Theorem~\ref{thm:EulerPhi} implies that a necessary condition for
this equality to hold is that
\[ E^{k} \phi_{m} - E^{m} \phi_{k} = (m-k) \phi_{k+m-1}. \]
In this section, we prove that this condition holds for all manifolds,
i.e. Theorem~\ref{thm:phiVir} is true.

We begin with the following
\begin{lem} \label{lem:derEulerH}
Let $H = \sum_{\beta} \gamma_{\beta} \bullet \gamma^{\beta}$.
For any $\alpha$, we have
\begin{eqnarray*}
 E^{2} \left<\left< \gamma_{\alpha} E^{k} H \right>\right>_{0, s}
 & = & (b_{\alpha}- b_{1}+k) 
	\left<\left< \gamma_{\alpha} E^{k+1} H \right>\right>_{0, s} \\
 &&   + \sum_{\mu} b_{\mu} \left<\left< \gamma_{\alpha} E \gamma^{\mu}
		 \right>\right>_{0, s}
	\left<\left< \gamma_{\mu} E^{k} H \right>\right>_{0, s} \\
 &&   - \sum_{\mu} b_{\mu} x_{1}^{\mu} 
	\left<\left< \gamma_{\mu} \left(\gamma_{\alpha} \bullet E^{k}
		\right) H \right>\right>_{0, s}.
\end{eqnarray*}
\end{lem} \label{lem:E2H}
{\bf Proof}: Since
\begin{eqnarray*}
 E^{2} \left<\left< \gamma_{\alpha} E^{k} H \right>\right>_{0, s}
&=& \sum_{\sigma, \beta} \left\{ E^{2} 
  \left<\left< \gamma_{\alpha} E^{k} \gamma^{\sigma} \right>\right>_{0, s}
	\right\}
  \left<\left< \gamma_{\sigma} \gamma_{\beta} \gamma^{\beta}
		 \right>\right>_{0, s} \\
&& + \sum_{\sigma, \beta} 
  \left<\left< \gamma_{\alpha} E^{k} \gamma^{\sigma} \right>\right>_{0, s}
  \left\{ E^{2} \left<\left< \gamma_{\sigma} \gamma_{\beta}
		\gamma^{\beta} \right>\right>_{0, s}
	\right\},
\end{eqnarray*}
the lemma follows by 
applying formula (\ref{eqn:derPowerEuler}) to the first term and
Lemma~\ref{lem:4ptEk} to the second term and then simplifying the resulting 
expression.
$\Box$

We also need the following
\begin{lem} \label{lem:E2prod}
\begin{eqnarray*}
&&  E^{2} \left\{ \sum_{i=0}^{k-1} \sum_{\alpha, \beta} b_{\alpha} b_{\beta} 
	\left<\left<  \gamma_{\alpha} E^{i} \gamma^{\beta}
		\right>\right>_{0, s}
	\left<\left<  \gamma_{\beta} E^{k-1-i} \gamma^{\alpha}
		\right>\right>_{0, s} \right\} \\
& = & k \sum_{\beta} b_{\beta}^{2} 
	\left<\left<  E^{k} \gamma_{\beta}  \gamma^{\beta}
		\right>\right>_{0, s}
  + (k-2)\sum_{i=0}^{k-1} \sum_{\alpha, \beta} b_{\alpha} b_{\beta} 
	\left<\left<  \gamma_{\alpha} E^{i+1} \gamma^{\beta}
		\right>\right>_{0, s}
	\left<\left<  \gamma_{\beta} E^{k-1-i} \gamma^{\alpha}
		\right>\right>_{0, s}. 
\end{eqnarray*}
\end{lem}
{\bf Proof}:
First observe that
\begin{eqnarray*}
&&  E^{2} \left\{ \sum_{i=0}^{k-1} \sum_{\alpha, \beta} b_{\alpha} b_{\beta} 
	\left<\left<  \gamma_{\alpha} E^{i} \gamma^{\beta}
		\right>\right>_{0, s}
	\left<\left<  \gamma_{\beta} E^{k-1-i} \gamma^{\alpha}
		\right>\right>_{0, s} \right\} \\
& = & 2\sum_{i=0}^{k-1} \sum_{\alpha, \beta} b_{\alpha} b_{\beta} 
	\left\{ E^{2} \left<\left<  \gamma_{\alpha} E^{i} \gamma^{\beta}
		\right>\right>_{0, s} \right\}
	\left<\left<  \gamma_{\beta} E^{k-1-i} \gamma^{\alpha}
		\right>\right>_{0, s}. 
\end{eqnarray*}
After applying formula (\ref{eqn:derPowerEuler}), we can simplify the
expression by using identities
\begin{eqnarray*}
&&   \sum_{i=0}^{k-1} \sum_{\alpha, \beta, \mu} b_{\alpha} b_{\beta} b_{\mu}
	\left<\left<  \gamma_{\alpha} E \gamma^{\mu}
		\right>\right>_{0, s}
	\left<\left<  \gamma_{\mu} E^{i} \gamma^{\beta}
		\right>\right>_{0, s}
	\left<\left<  \gamma_{\beta} E^{k-1-i} \gamma^{\alpha}
		\right>\right>_{0, s}  \\
& = &  \sum_{i=0}^{k-1} \sum_{\alpha, \beta, \mu} b_{\alpha} b_{\beta} b_{\mu}
	\left<\left<  \gamma_{\alpha} E^{i} \gamma^{\mu}
		\right>\right>_{0, s}
	\left<\left<  \gamma_{\mu} E \gamma^{\beta}
		\right>\right>_{0, s}
	\left<\left<  \gamma_{\beta} E^{k-1-i} \gamma^{\alpha}
		\right>\right>_{0, s},
\end{eqnarray*}
and
\begin{eqnarray*}
&&  \sum_{i=0}^{k-1} \sum_{\alpha, \beta} b_{\alpha}^{2} b_{\beta} 
	\left<\left<  \gamma_{\alpha} E^{i+1} \gamma^{\beta}
		\right>\right>_{0, s}
	\left<\left<  \gamma_{\beta} E^{k-1-i} \gamma^{\alpha}
		\right>\right>_{0, s}  \\
& = &  \sum_{i=0}^{k-1} \sum_{\alpha, \beta} b_{\alpha} b_{\beta}^{2} 
	\left<\left<  \gamma_{\alpha} E^{i+1} \gamma^{\beta}
		\right>\right>_{0, s}
	\left<\left<  \gamma_{\beta} E^{k-1-i} \gamma^{\alpha}
		\right>\right>_{0, s}.
\end{eqnarray*}
These two identities are obtained by interchanging indices.
We thus obtain
\begin{eqnarray*}
&&  E^{2} \left\{ \sum_{i=0}^{k-1} \sum_{\alpha, \beta} b_{\alpha} b_{\beta} 
	\left<\left<  \gamma_{\alpha} E^{i} \gamma^{\beta}
		\right>\right>_{0, s}
	\left<\left<  \gamma_{\beta} E^{k-1-i} \gamma^{\alpha}
		\right>\right>_{0, s} \right\} \\
& = & 2\sum_{i=0}^{k-1} \sum_{\alpha, \beta} i b_{\alpha} b_{\beta} 
	\left<\left<  \gamma_{\alpha} E^{i+1} \gamma^{\beta}
		\right>\right>_{0, s}
	\left<\left<  \gamma_{\beta} E^{k-1-i} \gamma^{\alpha}
		\right>\right>_{0, s}. 
\end{eqnarray*}
The lemma then follows from the identity
\begin{eqnarray*}
&& \sum_{i=0}^{k-1} \sum_{\alpha, \beta} i b_{\alpha} b_{\beta} 
	\left<\left<  \gamma_{\alpha} E^{i+1} \gamma^{\beta}
		\right>\right>_{0, s}
	\left<\left<  \gamma_{\beta} E^{k-1-i} \gamma^{\alpha}
		\right>\right>_{0, s} \\
& = & k \sum_{\beta} b_{\beta}^{2} 
	\left<\left<  E^{k} \gamma_{\beta}  \gamma^{\beta}
		\right>\right>_{0, s}
 +\sum_{j=0}^{k-1} \sum_{\alpha, \beta}(k-2-j) b_{\alpha} b_{\beta} 
	\left<\left<  \gamma_{\alpha} E^{j+1} \gamma^{\beta}
		\right>\right>_{0, s}
	\left<\left<  \gamma_{\beta} E^{k-1-j} \gamma^{\alpha}
		\right>\right>_{0, s}. 
\end{eqnarray*}
This identity is obtained by substituting $k-2-j$ for $i$ and
interchanging indices.
$\Box$

We can now prove a special case of Theorem~\ref{thm:phiVir}.
\begin{pro} \label{pro:phi2Vir}
\[ E^{k} \phi_{2} - E^{2} \phi_{k} = (2-k) \phi_{k+1}.\]
\end{pro}
{\bf Proof}:
By Lemma~\ref{lem:sEulerCorr}, 
\[ \phi_{2} = - \frac{1}{24} \sum_{\alpha, \beta} x_{1}^{\alpha}
	(b_{\alpha} - b_{1}) 
	\left<\left< \gamma_{\alpha} \gamma_{\beta} \gamma^{\beta}
			  \right>\right>_{0, s}
	+ \frac{1}{2} \sum_{\beta} \left\{ b_{\beta} (1-b_{\beta})
		- \frac{b_{1}+1}{6} \right\}
	\left<\left< \gamma_{\beta} \gamma^{\beta} \right>\right>_{0, s}.
\]
Using formula (\ref{eqn:derCorr}) and Lemma~\ref{lem:4ptEk}, we
can express $E^{k} \phi_{2}$ in terms of products of 3-point
functions. On the other hand, using Theorem~\ref{thm:phik}, formula
(\ref{eqn:derPowerEuler}), Lemma~\ref{lem:derEulerH} and
Lemma~\ref{lem:E2prod}, we can express $E^{2} \phi_{k}$
in terms of products of 3-point functions. Combine the two
expressions together and simplifying it, we obtain the desired
formula.
$\Box$
 
Now we are ready to prove Theorem~\ref{thm:phiVir}.

\noindent
{\bf Proof of Theorem~\ref{thm:phiVir}}:
We prove this theorem by induction on $\min \{m, \, k\}$. 
Without loss of generality, we may assume that $m \leq k$.

If $m = 0$, equation (\ref{eqn:phiVir}) is equivalent to
$\gamma_{1} \phi_{k} = k \phi_{k-1}$. This equality holds
trivially when $k = 0$ or $k=1$. When $k = 2$, it follows
from formula (\ref{eqn:derCorr}), Lemma~\ref{lem:StringCorr},
and the following formula (cf. \cite{Bor})
\[ \frac{1}{2} \sum_{\beta} b_{\beta} (1 - b_{\beta})
 - \frac{b_{1} + 1}{12} \chi (V) 
	= - \frac{1}{12} \int_{V} c_{1}(V) \cup c_{d-1}(V).\]
Note that this is the reason why $b_{\alpha}$ is defined
in terms of the holomorphic dimension of $\gamma_{\alpha}$
rather than a half of the real dimension of $\gamma_{\alpha}$
as proposed in \cite{EHX2}. For $k > 2$, the equality follows from
Theorem~\ref{thm:phik}, formula (\ref{eqn:derCorr}), the fact that
$\nabla_{\gamma_{1}} E^{k} = [\gamma_{1}, \, E^{k}] = k E^{k-1}$,
and Lemma~\ref{lem:StringCorr}. 

Assume that equality (\ref{eqn:phiVir}) holds for $m \leq n$. We
want to show that it also holds for $m = n+1$. In fact for any 
$k$, by equation (\ref{eqn:VirEuler}) and 
Proposition~\ref{pro:phi2Vir}, we have
\[
E^{n+1} \phi_{k} - E^{k} \phi_{n+1}
=  \frac{1}{n-2} \left\{ \left(E^{2}E^{n} - E^{n}E^{2}\right) \phi_{k} - 
		E^{k} \left(E^{2} \phi_{n} - E^{n} \phi_{2} \right) \right\}.
\]
By the induction hypothesis, 
$E^{n} \phi_{k} = E^{k} \phi_{n} + (k-n) \phi_{n+k-1}$,
and by Proposition~\ref{pro:phi2Vir},
$E^{2} \phi_{k} = E^{k} \phi_{2} + (k-2) \phi_{k+1}$.
Therefore, by equation (\ref{eqn:VirEuler}), we have
\[
E^{n+1} \phi_{k} - E^{k} \phi_{n+1}
=  \frac{1}{n-2} \left\{ (k-2) \left(E^{k+1} \phi_{n} 
						- E^{n} \phi_{k+1}	\right)  
		+(k-n) \left(E^{2} \phi_{n+k-1} - E^{n+k-1} \phi_{2} \right) 
	\right\}.
\]
Using the induction hypothesis and Proposition~\ref{pro:phi2Vir} again,
we have 
\[ E^{n+1} \phi_{k} - E^{k} \phi_{n+1} = (k-n-1) \phi_{n+k}.\]
This proves the theorem.
$\Box$

We can use Theorem~\ref{thm:phiVir} to construct a representation
of the Lie algebra spanned by $\{ E^{k} \mid k \geq 0\}$ in the
following way. 
Let 
\begin{equation} \label{eqn:hdef}
h_{k} := \left<\left< E^{k} \right>\right>_{1, s} - \phi_{k}.
\end{equation}
By Theorem~\ref{thm:g0g1G0} and the definition of $\phi_{k}$, 
$h_{0} = h_{1} = 0$ and
\begin{equation} \label{eqn:Ekh}
h_{k} = \frac{k}{2} E^{k-1} h_{2}.
\end{equation}
More generally, we have the following
\begin{lem} \label{lem:hVirRep} For all $k \geq 0$ and $m > 0$, 
\[ E^{k} \, \, \frac{h_{m}}{m} = (m-1) \, \,  \frac{h_{m+k-1}}{m+k-1}.\]
\end{lem}
{\bf Proof}: Theorem~\ref{thm:phiVir} and formula (\ref{eqn:VirEuler})
imply \[ E^{k} h_{m} - E^{m} h_{k} = (m-k) h_{k+m-1} \] for all $m$ and
$k$. Using this formula, one can show that the equation
\[ E^{k} \, \, \frac{h_{m}}{m} = (m-1) \, \,  \frac{h_{m+k-1}}{m+k-1}\]
is equivalent to the equation
\[ E^{m} \, \, \frac{h_{k}}{k} = (k-1) \, \,  \frac{h_{m+k-1}}{m+k-1}.\]
Formula (\ref{eqn:Ekh}) says that the lemma is true if $\min \{m, \, k\} = 2$.
By formulas (\ref{eqn:Ekh}) and (\ref{eqn:VirEuler}), we have
\[ E^{k} h_{m} = \frac{m}{2} \left\{ E^{m-1} \frac{2}{k+1} h_{k+1}
	+ (m-k-1) \frac{2}{m+k-1} h_{m+k-1} \right\}. \] 
The lemma then follows from  induction on $\min \{m, \, k\}$.
$\Box$

Lemma~\ref{lem:hVirRep} tells us that the linear span of
$\{ h_{k} \mid k \geq 2 \}$ gives a representation of the Lie algebra
spanned by $\{ E^{k} \mid k \geq 0\}$. Theorem~\ref{thm:phi2g1Vir} means that
the genus-1 Virasoro conjecture holds if and only if $h_{2} = 0$, which is
equivalent to say that this representation is trivial.

\section{Some sufficient conditions for genus-1 Virasoro conjecture}
\label{sec:SuffVir}

In an open subset of the small phase space, $\{ E^{k} \mid k \geq 0 \}$
defines an integrable distribution. Each leaf of this distribution is
a finite dimensional manifold. Fix one leaf of this distribution.
There exists an integer $n$ such that $\{ E^{k} \mid 0 \leq k \leq n \}$
are linearly independent and there are
smooth functions
$f_{i}$, $0 \leq i \leq n$, on the leaf such that 
\[ E^{n+1} = \sum_{i=0}^{n} f_{i} E^{i}. \]
Since $E^{k+n+1} = E^{k} \bullet E^{n+1}$, 
we have
\begin{equation} \label{eqn:wrapEuler2}
 E^{k+n+1} = \sum_{i=0}^{n} f_{i} E^{k+i}, 
\end{equation}
for every $k \geq 0$.
For later applications,  we need to compute 
$E^{k} f_{i}$.  
\begin{lem} \label{lem:fDer}
{\rm (i)} $E^{0} f_{i} = - (i+1) f_{i+1}$ {\rm for} $ 0 \leq i \leq n-1$, 
	{\rm and}
	$E^{0} f_{n} = n+1$. 

\hspace{48pt} {\rm (ii)} $E f_{i} = (n+1-i) f_{i}$ {\rm for} 
	$0 \leq i \leq n$ . 

\hspace{45pt}
{\rm (iii)} $ E^{2} f_{0} = f_{n} f_{0}$ {\rm and} 
	$E^{2} f_{i} = (n-i+2)f_{i-1} + f_{n} f_{i}$ {\rm for}
	 $1 \leq i \leq n$.

\hspace{48pt} {\rm (iv)} For $k > 0$, 
\[ \hspace{80pt} \left\{ \begin{array}{l}
		E^{k} f_{0} = f_{0} E^{k-1}f_{n}, \\ 
	E^{k} f_{i} = f_{i} E^{k-1}f_{n} + E^{k-1}f_{i-1} \, \, \,
	{\rm for} \, \, \, 1 \leq i \leq n.
		\end{array} \right.
\]
\end{lem}
\begin{rem} \label{rem:fDer}
{\rm 
Lemma~\ref{lem:fDer} (i) and (iv)
tell us that at each point, $E^{m} f_{j}$ is completely determined by the
values of $f_{0}, \ldots, f_{n}$ at that point. }
\end{rem}
{\bf Proof of Lemma~\ref{lem:fDer}}: 
We first prove formula (iv). By formulas (\ref{eqn:VirEuler}) and
(\ref{eqn:wrapEuler2}), for $0 \leq m \leq k$, 
\begin{eqnarray}
(n+1+2m-k) E^{n+k} & = & \left[ E^{k-m}, \, \, E^{n+m+1} \right]
		 = \left[ E^{k-m}, \, \, \sum_{i=0}^{n} f_{i} E^{i+m} \right] 
			\nonumber \\
	&= & \sum_{i=0}^{n} \left( E^{k-m} f_{i} \right) E^{i+m} 
		 + \sum_{i=0}^{n} f_{i} \left[ E^{k-m},  E^{i+m} \right] 
			\nonumber \\
	&=& \sum_{i=0}^{n} \left( E^{k-m} f_{i} \right) E^{i+m} 
		+  \sum_{i=0}^{n} f_{i} (i+2m-k) E^{i+k-1}. \label{eqn:fDerPf}
\end{eqnarray}
Using the fact that 
$2m E^{k+n} = \sum_{i=0}^{n} 2m f_{i} E^{i+k-1}$, we obtain
\[
	\sum_{i=0}^{n} \left( E^{k-m} f_{i} \right) E^{i+m} 
	= (n+1-k) E^{n+k} -   \sum_{i=0}^{n} f_{i} (i-k) E^{i+k-1}. 
\]
Since the right hand side of this equation does not depend on $m$,
so does the left hand side. Therefore we have
\begin{equation} \label{eqn:RecZk}
 \sum_{i=0}^{n} \left( E^{k} f_{i} \right) E^{i} 
	= \sum_{i=0}^{n} \left( E^{k-m} f_{i} \right) E^{i+m},
\end{equation}
for all $0 \leq m \leq k$. 
In the special case $m=1$, we have
\[ \sum_{i=0}^{n} \left( E^{k} f_{i} \right) E^{i} 
	= \sum_{i=0}^{n} \left( E^{k-1} f_{i} \right) E^{i+1}. \]
Replacing $E^{n+1}$ on the right hand side of this equality by
$ \sum_{i=0}^{n} f_{i} E^{i} $ and using the fact that
$\{ E^{0}, \ldots, E^{n} \}$ are linearly independent, we obtain
formula (iv). 

Formula (i) is obtained from (\ref{eqn:fDerPf}) by setting $k = m = 0$.
Formula (ii) and (iii) are obtained by using (i) and the recursion
formula (iv).
$\Box$ 

Now we come back to the Virasoro conjecture. 
As pointed out in the introduction, a necessary condition for the 
genus-1 Virasoro conjecture to hold is the validity of formula
(\ref{eqn:philinear}), i.e.
\[  \phi_{n+1} = \sum_{k=0}^{n} f_{k} \phi_{k}. \]
This condition implies the
following
\begin{lem} \label{lem:philinn+k}
If formula (\ref{eqn:philinear}) is correct, then
\[  \phi_{m+n+1} = \sum_{k=0}^{n} f_{k} \phi_{m+k}, \]
for all $m \geq 0$.
\end{lem}
{\bf Proof}:
By formula (\ref{eqn:phiVir}),
\begin{eqnarray*}
&& E^{2} \left\{ \phi_{m+n+1} - \sum_{j=0}^{n} f_{j} \phi_{m+j} \right\} \\
& = & E^{m+n+1} \phi_{2} + (m+n-1) \phi_{m+n+2}  \\
&&	 - \sum_{j=0}^{n} \left( E^{2}f_{j} \right) \phi_{m+j} 
 	- \sum_{j=0}^{n} f_{j} \left\{E^{m+j}\phi_{2} + (m+j-2)\phi_{m+j+1}
		 \right\}. 
\end{eqnarray*}
Using Lemma~\ref{lem:fDer} (iii) and the formula (\ref{eqn:wrapEuler2}),
we obtain
\begin{eqnarray*}
&&  E^{2} \left\{ \phi_{m+n+1} - \sum_{j=0}^{n} f_{j} \phi_{m+j} \right\} \\
&=& (m+n-1) \left\{ \phi_{m+n+2} - 
		\sum_{j=0}^{n} f_{j} \phi_{m+j+1} \right\}
	  + f_{n} \left\{ \phi_{m+n+1} - 
		\sum_{j=0}^{n} f_{j} \phi_{m+j} \right\}. 
\end{eqnarray*}
The lemma then follows by induction on $m$. 
$\Box$

By Theorem~\ref{thm:phi2g1Vir}, to prove the genus-1 Virasoro conjecture
we only need to show that
$h_{2} := \left<\left<  E^{2} \right>\right>_{1, s} - \phi_{2} = 0$.
We first prove the following
\begin{pro} \label{pro:Yh0}
Let 
\[ Z_{k} := \sum_{i=0}^{n} \left( E^{k} f_{i} \right) E^{i}. \]
If equality (\ref{eqn:philinear}) holds, then
$Z_{k} h_{2} = 0 $ for all $k \geq 0$.
\end{pro}
\begin{rem} \label{rem:ZkZ0}
{\rm Formula~(\ref{eqn:RecZk}) implies that $Z_{k} = E^{k} \bullet Z_{0}$.}
\end{rem}
{\bf Proof of Proposition~\ref{pro:Yh0}}:
Setting $m = k$ in formula~(\ref{eqn:fDerPf})
 and using formula~(\ref{eqn:RecZk}),
we obtain
\[ Z_{k} = (n+k+1) E^{n+k} -  \sum_{i=0}^{n} (i+k) f_{i}  E^{i+k-1}. \]
Therefore by Lemma~\ref{lem:hVirRep},
\[ Z_{k} h_{2} = 2 h_{n+k+1} - \sum_{i=0}^{n} 2 f_{i}  h_{i+k}. \]
The right hand side of this equality is  equal to $0$ because of
formula (\ref{eqn:wrapEuler2}) and Lemma~\ref{lem:philinn+k}.
$\Box$

An immediate consequence of this proposition is Theorem~\ref{thm:nondegVir}.

{\bf Proof of Theorem~\ref{thm:nondegVir}}: By Lemma~\ref{lem:hVirRep},
$E^{0} h_{2} = 0$. If for some positive integer $m$, 
 $E^{m}$ is contained in the span of
$\{ E^{0}, Z_{k} \mid k \geq 0 \}$, then by Proposition~\ref{pro:Yh0},
$E^{m} h_{2} = 0$. By Lemma~\ref{lem:hVirRep}, 
\[ h_{m+1} = \frac{m+1}{2} \, \,  E^{m}h_{2} = 0.\] 
If $m \geq 1$, then 
repeatedly taking derivatives (by $m-1$ times) of $h_{m+1}$ along the 
direction $E^{0}$ and
using Lemma~\ref{lem:hVirRep}, we obtain that $h_{2} = 0$.  The
theorem then follows from Theorem~\ref{thm:phi2g1Vir}.
$\Box$

To apply Theorem~\ref{thm:nondegVir}, we need to know which manifolds have
non-degenerate quantum cohomology. In the rest of this paper, we
discuss some sufficient conditions for the non-degeneracy of the
quantum cohomology. 
To this end,  it is interesting to know how large is the vector space spanned
by $\{Z_{k} \mid k \geq 0\}$. We first notice that by Remark~\ref{rem:ZkZ0}
and formula~(\ref{eqn:wrapEuler2}),
\[ Z_{n+1+k} = \sum_{i=0}^{n} f_{i} Z_{i+k}. \]
Therefore, at each point, $\{Z_{k} \mid k \geq 0\}$ and
$\{ Z_{k} \mid 0 \leq k \leq n \}$ span the same vector space.
The following lemma gives us a sense on how large this vector
space might be.
\begin{lem} \label{lem:NonDegSS}
At each point $t$, 
\[ {\rm span}\{Z_{k}(t) \mid 0 \leq k \leq n\}
= {\rm span}\{E^{k}(t) \mid 0 \leq k \leq n\} \]
if and only if the polynomial in $x$
\[ p_{t}(x) = x^{n+1} - \sum_{i=0}^{n} f_{i}(t) x^{i} \]
has no multiple roots.
\end{lem}
{\bf Proof}:
The derivative of $p_{t}(x)$ with respect to $x$ is
\[ p_{t}^{'}(x) = (n+1)x^{n} - \sum_{i=0}^{n-1} (i+1)f_{i+1}(t) x^{i}. \]
Note that the coefficients of $p_{t}^{'}(x)$ are the same as
the coefficients of $Z_{0}(t)$. The resultant of polynomials
$p_{t}(x)$ and $p_{t}^{'}(x)$ is the determinant of the following
$(2n+1) \times (2n+1)$ matrix
\[ \footnotesize{ \left( \begin{array}{ccccccccc} 
		1, & -f_{n}, & -f_{n-1}, & \cdots \cdots \cdots,
			 & -f_{1} & -f_{0}, & & & \\
		 & 1, & -f_{n}, & -f_{n-1}, & \cdots \cdots,
			 & -f_{1}& -f_{0}, & &  \\
		& & \cdots \cdots \cdots \cdots & \cdots \cdots \cdots \cdots
			 & \cdots \cdots & 
			\cdots \cdots& \cdots \cdots \cdots \cdots& 
			\cdots \cdots  & \\
		&&  & 1, & -f_{n}, & -f_{n-1}, & \cdots \cdots \cdots, 
				&-f_{1} & -f_{0}  \\
		&&&&&&&& \\
		n+1, & -nf_{n}, & -(n-1)f_{n-1}, & \cdots \cdots \cdots, 
				& -f_{1}, & & & &\\
		 & n+1, & -nf_{n}, & -(n-1)f_{n-1}, & \cdots \cdots, & -f_{1}, 
			& & & \\
		& & \cdots \cdots \cdots \cdots& 
			\cdots \cdots \cdots \cdots & 
			\cdots \cdots & \cdots \cdots & 
			\cdots \cdots \cdots \cdots & & \\
		& & & \cdots \cdots \cdots \cdots & \cdots \cdots & 
		\cdots \cdots & \cdots \cdots \cdots \cdots & \cdots \cdots & \\
		& &&  & n+1, & -nf_{n}, & -(n-1)f_{n-1}, & \cdots 
			\cdots, & -f_{1} 
	  \end{array} \right)_{,} }
\]
where non-zero entries of the first
$n$ rows are coefficients of $p_{t}(x)$ and non-zero 
entries of the last
$n+1$ rows are coefficients of $p_{t}^{'}(x)$.  Performing elementary
row transformations, we can transform this matrix to the following form
\[ \left( \begin{array} {cc}
		B & C \\ 0 & A \end{array} \right), \]
where $B$ is an $n \times n$ upper triangular matrix whose diagonal
entries are 1, and $A = (a_{i, j})$, $0 \leq i, j \leq n$,
 is an $(n+1) \times (n+1)$ matrix
whose entries are given by the  recursion formula
\[ \begin{array}{l}
	a_{n, 0} = n+1, \, \, \, a_{n, j} = - (n-j+1) f_{n-j+1} 
	\, \, \, {\rm for} \, \, \, 1 \leq j \leq n; \\ 
	{\rm For} \, \, \, 1 \leq i \leq n, \\
	 \hspace{20pt} \left\{ \begin{array}{l}
		a_{n-i, n} = f_{0} a_{n-i+1, 0}, \\
		a_{n-i, j} = f_{n-j} a_{n-i+1, 0} + a_{n-i+1, j+1},
		\, \, \, {\rm for} \, \, \, 0 \leq j \leq n-1.
		\end{array} \right.
   \end{array}
\]
Comparing this recursion formula with the recursion formula in 
Lemma~\ref{lem:fDer}, we obtain that $a_{i, j} = E^{n-i} f_{n-j}$
for all $i$ and $j$. Therefore $A$ is the coefficient matrix of
representing $\{Z_{n}, Z_{n-1}, \ldots, Z_{0} \}$ in terms of
$\{ E^{n}, E^{n-1}, \ldots, E^{0} \}$. Since the determinant of
$A$ is equal to the resultant of $p_{t}(x)$ and $p_{t}^{'}(x)$,
$A$ is invertible if and only if $p_{t}(x)$ has no multiple
roots. This proves the lemma.
$\Box$

Recall that a manifold has non-degenerate quantum cohomology if
there exists one $m > 0$ such that at generic points, 
 $E^{m}$ is contained in the span of $\{ E^{0}, Z_{0}, \ldots, Z_{n} \}$.
Observe that if the first $n$ columns of the matrix $A$ in the proof
of Lemma~\ref{lem:NonDegSS} has rank $n$, than $E^{m}$ is contained
in the span of $\{ E^{0}, Z_{0}, \ldots, Z_{n} \}$ for all $m \geq 0$.
Therefore such manifolds have non-degenerate quantum cohomology.
However, to compare non-degeneracy with semisimplicity, 
we only need the following
weaker result which corresponds to the case where $A$ has rank $n+1$. 
\begin{cor} \label{cor:NonDegSS}
If at generic points of the small phase space of a manifold $V$,  
the polynomial 
\[ p_{t}(x) = x^{n+1} - \sum_{i=0}^{n} f_{i}(t) x^{i} \]
has no multiple roots, then the quantum cohomology of $V$ is non-degenerate.
\end{cor}

In the case that the quantum cohomology of $V$ is semisimple, 
at generic points of the small phase space, $\{E^{k} \mid 0 \leq k \leq n\}$
form a basis of the tangent space of the small phase space. With respect
to this basis, the quantum multiplication by $E$ has the following
matrix representation
\[ \left( \begin{array}{cccccc}
	0 & 0 & 0 & \cdots & 0 & f_{0} \\
	1 & 0 & 0 & \cdots & 0 & f_{1} \\
	0 & 1 & 0 & \cdots & 0 & f_{2} \\
	0 & 0 & 1 & \cdots & 0 & f_{3} \\
	\cdots & \cdots &\cdots & \cdots & \cdots & \cdots \\
	0 & 0 & 0 & \cdots & 1 & f_{n}
	\end{array} \right). \]
The polynomial $p_{t}(x)$ in Corollary~\ref{cor:NonDegSS} is precisely
the characteristic polynomial of this matrix, and therefore
 has no multiple roots at semisimple points. Hence we have
\begin{cor} \label{cor:NonDegSS2}
If the quantum cohomology of a manifold $V$ is semisimple, then it is
also non-degenerate.
\end{cor}

Another sufficient condition for the non-degeneracy is the following
\begin{lem} \label{lem:nondegdim1}
If at every point of the small phase space, 
the dimension of the
vector space spanned by $\{ E^{k} \mid k \geq 0 \}$ 
is less than or equal
to 2, then the quantum cohomology is non-degenerate.
\end{lem}
{\bf Proof}:
The case where the dimension of the
vector space spanned by $\{ E^{k} \mid k \geq 0 \}$ is 1 is trivial
since $E^{k}$ is proportional to $E^{0}$ for every $k > 1$.
If the dimension of the
vector space spanned by $\{ E^{k} \mid k \geq 0 \}$ is 2, then
$E^{2} = f_{0}E^{0} + f_{1}E$ with $E^{0} f_{1} = 2$ 
(c.f. Lemma~\ref{lem:fDer} (i)). Hence 
\[ E = \frac{1}{2} \left\{ Z_{0} - \left(E^{0}f_{0}\right) E^{0} \right\}. \]
By definition, the quantum cohomology is non-degenerate.
$\Box$

Now we give some examples where the quantum cohomologies are
non-degenerate but not semisimple.
\begin{exa} \label{exa:curve}
{\rm Let $C_{g}$ be a complex curve of genus $g$.
We first consider only even dimensional cohomology classes 
$H^{\rm even}(C_{g}; {\Bbb C})$. Since 
 the dimension of the small phase space is 2, the quantum cohomology
is non-degenerate. However except for
$g = 0$, the quantum cohomology of $C_{g}$ is not semisimple since
its first Chern class is either zero or negative. In the case of complex
one dimensional tori, 
the Euler vector field is proportional to the identity element. Equality 
(\ref{eqn:philinear}) is trivially satisfied since $\phi_{1} = \phi_{0} = 0$.
However we can not apply Theorem~\ref{thm:nondegVir} yet because 
$C_{g}$ has non-trivial odd dimensional cohomology classes. The reason is
that if we do not consider all cohomology classes, we can not get
the Euler characteristic number in Borisov's formula (see the proof
of Theorem~\ref{thm:phiVir}). So Theorem~\ref{thm:phiVir}, which is a
 necessary condition for genus-1 Virasoro conjecture, 
does not hold if we only consider even dimensional cohomology classes.
 
To prove the genus-1 Virasoro conjecture for $C_{g}$ with $g > 0$, 
we have to consider the space of all cohomology classes 
$H^{*}(C_{g}; {\Bbb C})$. All theorems stated in the introduction can
be extended to manifolds with non-trivial odd dimensional cohomologies
without any difficulty. Therefore these theorems can be applied to
$C_{g}$.   Since there is no non-constant holomorphic
maps from a rational curve to $C_{g}$, the quantum cohomology of $C_{g}$
is the same as the ordinary cohomology. It follows that
$E^{2} = -(t^{1})^{2}E^{0} + 2t^{1} E$. Therefore the quantum cohomology
is again non-degenerate but not semisimple. Since the only genus-0 non-zero
Gromov-Witten invariants for $C_{g}$ are 3-point degree-0 invariants
which can be computed via cup products, it is straightforward to check that
$\phi_{2} = \frac{g-1}{6} t^{1}$. Therefore equality (\ref{eqn:philinear})
holds. By theorem~\ref{thm:nondegVir},
 the genus-1 Virasoro conjecture holds for all complex
curves. To our knowledge, this result is not known before.
}
\end{exa}

\begin{exa} \label{exa:K3}
{\rm
Let $V$ be a K3 surface. Let $\gamma_{1}$ be the identity element of 
the cohomology ring of $V$ and $\gamma_{N}$ be a non-zero element of
$H^{4}(V)$. Since $c_{1}(V) = 0$,  
\[ E = t^{1}\gamma_{1} - t^{N} \gamma_{N}.\]
Because of the selection rule and the puncture equation, on the small
phase space, any $k$-point function
involving $\gamma_{N}$ is zero if $k \geq 4$ and the only non-zero 3-point
function involving $\gamma_{N}$ is 
$\left<\left< \gamma_{N} \gamma_{1} \gamma_{1} \right>\right>_{0, s}$. In
particular, $\gamma_{N} \bullet \gamma_{N} = 0$. Therefore the vector
space spanned by $\{ E^{k} \mid k \geq 0 \}$ is of dimension 2. Hence
the quantum cohomology of $V$ is non-degenerate. It is not semisimple
since $c_{1}(V) = 0$. Moreover 
$ \gamma_{1} \phi_{2} = E^{0} \phi_{2} = 2 \phi_{1} = 0$ by formula
(\ref{eqn:phiVir}), and 
\begin{eqnarray*}
 \gamma_{N} \phi_{2} & = & 
	- \frac{1}{24} \sum_{\alpha} 
         \left<\left< \gamma_{N} E E \gamma_{\alpha} \gamma^{\alpha} 
		\right>\right>_{0, s} 
	 + \frac{1}{12} \sum_{\alpha} 
         \left<\left< \gamma_{N} E \gamma_{\alpha} \gamma^{\alpha} 
		\right>\right>_{0, s} \\
        && + \frac{1}{2}\sum_{\alpha}  \left(b_{\alpha} (1-b_{\alpha})
                 - \frac{b_{1}+1}{6}\right)
        \left<\left< \gamma_{N}\gamma_{\alpha} \gamma^{\alpha} 
	\right>\right>_{0, s} \\
	& = & 0. 
\end{eqnarray*}
Therefore, by formula
(\ref{eqn:phiVir}) again, $\phi_{2} = E \phi_{2} = 0$. Consequently
equality (\ref{eqn:philinear}) holds trivially. Hence the
genus-1 Virasoro conjecture  holds for K3 surfaces.
 For Calabi-Yau manifolds with
complex dimension bigger than 2, the Virasoro conjecture holds
for dimension reasons (c.f. \cite{G2}). Therefore we know
that the genus-1 Virasoro conjecture holds for all Calabi-Yau manifolds.

After the first version of this paper had been posted on the web,
the author was informed by Jim Bryan that the Virasoro conjecture for
K3 surfaces follows from deformation invariance of the virtual moduli
cycle and triviality of moduli space of stable maps for generic K3
surfaces except for the genus-1 case where the virtual moduli cycle
of degree 0 maps is non-trivial. In the genus-1 case, one can use
intersection theory on $M_{1, n}$ to prove it. He also informed the
author that F. Zahariev found a combinatorial proof to the genus-1
Virasoro conjecture for the 1-dimensional tori. 
}
\end{exa}
 
\begin{exa} \label{exa:Pn}
{\rm
As for the semisimplicity, the non-degeneracy can also be defined for an 
abstract Frobenius manifold in the same way. We consider the 
Frobenius manifold $M_{n} := H^{*}({\Bbb C}P^{n})$ where the Frobenius
algebra structure is given by the ordinary cohomology ring structure
at every point of $M_{n}$ (c.f. Example 1.5 in \cite{Du}). 
It is not semisimple since it has nilpotent elements at each point.
Let $\gamma$ be a non-zero element of $H^{2}({\Bbb C}P^{n})$.
Then $\{ \gamma^{k} \mid 0 \leq k \leq n \}$ form a basis of $M_{n}$,
where $\gamma^{k} = \underbrace{\gamma \cup \cdots \cup \gamma}_{k}$.
We denote the corresponding coordinates by $t_{k}$, $0 \leq k \leq n$.
(This notation is different from our convention before where superscripts
were used instead of subscripts.)
The Euler vector field on $M_{n}$ is given by 
$E = \sum_{k=0}^{n} (1-k) t_{k} \gamma^{k}$.
It is straightforward to verify that for $n \leq 3$,
the dimension of the vector space spanned by $\{ E^{k} \mid k \geq 0 \}$
is less than or equal to 2. Therefore, in this case,
 $M_{n}$ is non-degenerate. For $n = 4$ or $5$,
$E^{3} = t_{0}^{3} E^{0} - 3t_{0}^{2}E + 3t_{0} E^{2}$. Therefore
$Z_{0} = 3t_{0}^{2} E^{0} - 6t_{0}E + 3 E^{2}$ and 
$Z_{k} = t_{0}^{k} Z_{0}$ for $k \geq 1$. Therefore $M_{4}$ and $M_{5}$
are degenerate. Notice that in this example,
the polynomial in Corollary~\ref{cor:NonDegSS} is of the form
$p_{t}(x) = (x - t_{0})^{3}$. In fact, it is not hard to show that
in general, 
if the dimension of the vector space spanned by $\{E^{k} \mid k \geq 0\}$
is equal to 3, then the Frobenius manifold is non-degenerate unless
$p_{t}(x) = (x - g(t))^{3}$ for some function $g$. For $M_{6}$,
we have $E^{4} = -t_{0}^{4} E^{0} + 4t_{0}^{3}E- 6t_{0}^{2}E^{2} 
		+ 4t_{0} E^{3}$. Therefore
$Z_{0} = -4t_{0}^{3}E^{0} + 12 t_{0}^{2} E - 12t_{0}E^{2} + 4E^{3}$ and 
$Z_{k} = t_{0}^{k} Z_{0}$ for $k \geq 1$. Therefore $M_{6}$ is also
degenerate. 
}
\end{exa}

%\left<\left<  \right>\right>_{0, s}
%%%$$$$%%%%%%%%%%%%%%%%%%%%%%%%%%%%%%%%%%%%%%%%%%%%%%%%%%%%%%%%%%%%%%%%%%%

%%%%%%%%%%%%%%%%%%%%%%%%%%%%%%%%%%%%%%%%%%%%%%%%%%%%%%%%%%%%%%%%%%%
\vspace{30pt}

\[ 
\begin{array}{lll} 
&&  New \, \, \, address: \\
{\rm Department \, \, \, of \, \, \, Mathematics} & \hspace{80pt}& 
		{\rm Department \, \, \, of \, \, \, Mathematics} \\ 
{\rm Massachusetts \, \, \, Institute \, \, \, of \, \, \, Technology} & & 
		{\rm University \, \, \, of \, \, \, Notre \, \, \, Dame} \\
{\rm Cambridge, \, \, \, MA \, \, \, 02139} & & 
		{\rm Notre \, \, \, Dame, \, \, \, IN \, \, \, 46556} \\
{\rm USA} & & 
		{\rm USA} \\
&& \\
{\rm E-mail \, \, \, address:} & & {\rm E-mail \, \, \, address:} \\
    xbliu@math.mit.edu && \it xliu3@nd.edu
\end{array}
\]

\end{document}